\documentclass{article}

\newcommand{\sfZF}{{\sf ZF}}

\usepackage{proof}
\usepackage{latexsym}
\usepackage{amsfonts}
\usepackage{amssymb}
\usepackage{cite}
\usepackage{color}

\newcommand{\blem}{\begin{lemma}}
\newcommand{\elem}{\end{lemma}}
\newcommand{\bth}{\begin{theorem}}
\newcommand{\ethm}{\end{theorem}}
\newcommand{\benu}{\begin{enumerate}}
\newcommand{\eenu}{\end{enumerate}}
\newcommand{\bdes}{\begin{description}}
\newcommand{\edes}{\end{description}}
\newcommand{\bdf}{\begin{definition}}
\newcommand{\edf}{\end{definition}}
\newcommand{\bcor}{\begin{cor}}
\newcommand{\ecor}{\end{cor}}
\newcommand{\bprp}{\begin{proposition}}
\newcommand{\eprp}{\end{proposition}}
\newcommand{\bmlem}{\begin{mlemma}}
\newcommand{\emlem}{\end{mlemma}}
\newcommand{\bclm}{\begin{claim}}
\newcommand{\eclm}{\end{claim}}
\newcommand{\bprf}{{\bf Proof}.\hspace{2mm}}

\newcommand{\eprf}{\hspace*{\fill} $\Box$}

\newcommand{\beqn}{\begin{equation}}
\newcommand{\eeqn}{\end{equation}}
\newcommand{\beqnarr}{\begin{eqnarray}}
\newcommand{\eeqnarr}{\end{eqnarray}}
\newcommand{\beqnarrs}{\begin{eqnarray*}}
\newcommand{\eeqnarrs}{\end{eqnarray*}}

\newcommand{\spand}{\,\&\,}

\newtheorem{theorem}{Theorem}[section]
\newtheorem{definition}[theorem]{Definition}
\newtheorem{proposition}[theorem]{Proposition}
\newtheorem{lemma}[theorem]{Lemma}
\newtheorem{cor}[theorem]{Corollary}
\newtheorem{mlemma}[theorem]{Main Lemma}
\newtheorem{claim}[theorem]{Claim}

\newcommand{\alp}{\alpha}

\newcommand{\veps}{\varepsilon}
\newcommand{\del}{\delta}
\newcommand{\Del}{\Delta}
\newcommand{\ome}{\omega}

\newcommand{\bet}{\beta}
\newcommand{\gam}{\gamma}
\newcommand{\Gam}{\Gamma}
\newcommand{\kap}{\kappa}
\newcommand{\sig}{\sigma}
\newcommand{\Sig}{\Sigma}
\newcommand{\tht}{\theta}
\newcommand{\Tht}{\Theta}
\newcommand{\lam}{\lambda}
\newcommand{\Lam}{\Lambda}
\newcommand{\vphi}{\varphi}

\newcommand{\fal}{\forall}
\newcommand{\exi}{\exists}

\newcommand{\Rarw }{\Rightarrow}

\newcommand{\lrarw}{\leftrightarrow}
\newcommand{\Lrarw}{\Leftrightarrow}

\newcommand{\calh}{{\cal H}}

\newcommand{\calL}{{\cal L}}
\newcommand{\calk}{{\cal K}}

\newcommand{\calP}{{\cal P}}

\newcommand{\rk}{\mbox{{\rm rk}}}

\newcommand{\la}{\langle}
\newcommand{\ra}{\rangle}

\newcommand{\setm}{\setminus}

\title{
Lifting proof theory to the countable ordinals II: second-order indescribable cardinals
}
\author{Toshiyasu Arai
\\
Graduate School of Science,
Chiba University
\\
1-33, Yayoi-cho, Inage-ku,
Chiba, 263-8522, JAPAN
\\
tosarai@faculty.chiba-u.jp
}
\date{}
\begin{document}
\maketitle

\begin{abstract}
We show that the existence of a $\Pi^{1}_{N}$-indescribable cardinal over the Zermelo-Fraenkel's set theory ${\sf ZF}$
is proof-theoretically reducible to iterations of 
Mostowski collapsings and 
lower Mahlo operations.
Furthermore we describe
a proof-theoretic bound on definable countable ordinals
whose existence is provable from the existence of second order indescribable cardinals over $\sfZF$.
\end{abstract}

\section{Introduction}\label{sect:introduction}

In \cite{liftupK} we showed that
the existence of a weakly compact, i.e., $\Pi^{1}_{1}$-indescribable cardinal over the Zermelo-Fraenkel's set theory ${\sf ZF}$
is proof-theoretically reducible to iterations of Mostowski collapsings and Mahlo operations,
while in \cite{liftupZF} we describe a proof-theoretic bound on definable countable ordinals in $\sfZF$.
In this paper we do the same reductions for the existence of a $\Pi^{1}_{N}$-indescribable cardinal, 
cf. Theorems \ref{thm:1k} and \ref{thm:2} in the last section \ref{sec:theorem}.

Let $ORD$ denote the class of all ordinals, $A\subset ORD$ and $\alp$ a limit ordinal.
$A$ {\rm is said to be} $\Pi^{1}_{n}${\it -indescribable\/} {\rm in} $\alp$ {\rm iff}
{\rm for any} $\Pi^{1}_{n}${\rm -formula} $\vphi(X)$ {\rm and any} $B\subset ORD$, {\rm if}
$\la L_{\alp},\in, B\cap\alp\ra\models\vphi(B\cap\alp)$,
{\rm then there exists a} $\bet\in A\cap\alp$ {\rm such that}
$\la L_{\bet},\in, B\cap\bet\ra\models\vphi(B\cap\bet)$.
Let us write
\[
\alp\in M_{n}(A):\Lrarw A \mbox{ is } \Pi^{1}_{n}\mbox{-indescribable in } \alp
.\]
Also $\alp$ {\rm is said to be} $\Pi^{1}_{n}${\it -indescribable\/} {\rm iff}
$\alp$ 
{\rm is} $\Pi^{1}_{n}${\rm -indescribable in} $\alp$.

It is not hard to extend the reduction in \cite{liftupK}
to $\Pi^{1}_{n}$-indescribability.
Namely over $\sfZF+(V=L)$ 
the existence of a $\Pi^{1}_{n}$-indescribable cardinal
is shown to be proof-theoretically reducible to iterations of 
Mostowski collapsings and 
the operation $M_{n-1}$,
since a similar reduction has been done for first-order reflecting ordinals
in \cite{consv}.

In this paper we aim a proof-theoretic reduction of $\Pi^{1}_{n}$-indescribability in terms of iterations of Mostowski collapsings and 
the operations $M_{i}$ for $i<n$.
Though such a reduction was done for recursive analogues, i.e., $\Pi_{n}$-reflecting ordinals and recursively Mahlo operations
in \cite{LMPS, WFnonmon2,PTPiN},
our approach is simpler.

First in \cite{liftupK} we rely on a result by R. Jensen\cite{Jensen}, which is the case $n=1$ in
a recent result (Theorem \ref{th:BagariaMagidorSakai} below) 
due to J. Bagaia, M. Magidor and H. Sakai\cite{Sakai}.

\bdf\label{df:BagariaMagidorSakai}
{\rm Let} $A\subset ORD$, {\rm and} $\alp$ {\rm a limit ordinal.}
\benu
\item
$A$ {\rm is said to be} $0$-stationary {\rm in} $\alp$ {\rm iff}
$\sup(A\cap\alp)=\alp$.

\item
{\rm For} $n>0$, $A$ {\rm is said to be} $n$-stationary {\rm in} $\alp$ {\rm iff}
{\rm for every} $m<n$ {\rm and every} $S\subset ORD$,
{\rm if} $S$ {\rm is} $m${\rm -stationary in} $\alp$, {\rm then there exists a} $\bet\in A\cap\alp$
{\rm such that} $S$ {\rm is} $m${\rm -stationary in} $\bet$.

\item
$\alp$ {\rm is said to be} $n$-stationary {\rm iff} $\alp=\{\bet\in ORD:\bet<\alp\}$ 
{\rm or} $ORD$ 
{\rm is $n$-stationary in} $\alp$.

\eenu
\edf
Note that $A$ is $1$-stationary in $\alp$ of uncountable cofinality iff $A\cap\alp$ is stationary in $\alp$,
i.e., $A$ meets every club subset of $\alp$.

\bth\label{th:BagariaMagidorSakai}{\rm (J. Bagaia, M. Magidor and H. Sakai\cite{Sakai})}
Let $\kap$ be a regular uncountable cardinal, and $A\subset ORD$.
For each $n>0$,
$A$ is $(n+1)$-stationary in $\kap$ 
iff $A$ is $\Pi^{1}_{n}$-indescribable in $\kap$,
over $\sfZF+(V=L)$.
\end{theorem}

Although the theorem is suggestive, we don't rely on it in this paper.

Second our classes $Mh_{k,n}(\vec{\alp})[\Tht]$ defined in Definition \ref{df:psivecalp} 
to resolve or approximate $\Pi^{1}_{N}$-indescribability are defined from 
finite \textit{sequences} of ordinals $\vec{\alp}$.
In \cite{LMPS, WFnonmon2,PTPiN}, our ramification process is akin to a tower, i.e., has an exponential structure.
Here we simplify the complicated process in terms of sequences.
Also cf. an ordinal analysis for first-order reflection using reflection configurations by Pohlers and Stegert \cite{Pohlers}.
\\

Let us mention the contents of this paper. 

In the next section \ref{sect:ordinalnotation} 
iterated Skolem hulls $\mathcal{H}_{\alpha,n}(X)$ of sets $X$ of ordinals, ordinals $\Psi_{\kappa,n}\gamma$ for regular ordinals 
$\kappa\,(\mathcal{K}<\kappa\leq I)$,
and classes $Mh_{k,n}(\vec{\alpha})[\Theta]$ are defined for finite sequences $\vec{\alp}$ of ordinals and
finite sets $\Theta$ of ordinals.
It is shown that for each $k<N$ and each $n,m<\omega$, $(\mathcal{K} \mbox{ {\rm is a }} \Pi^{1}_{N}\mbox{{-indescribable cardinal}}) \to \mathcal{K}\in Mh_{k,n}((\omega_{m}(I+1),\ldots,\omega_{m}(I+1)))[\emptyset]$ in
${\sf ZF}+(V=L)$.
In the third section \ref{sect:Ztheory} we introduce a theory for $\Pi^{1}_{N}$-indescribable cardinals,
which are equivalent to ${\sf ZF}+(V=L)+(\mathcal{K} \mbox{ {\rm is a }} \Pi^{1}_{N}\mbox{-indescribable cardinal}) $.

In the section \ref{sect:controlledOme}
cut inferences are eliminated from operator controlled derivations of $\Pi^{1}_{k}$-sentences $\varphi^{V_{\mathcal{K}}}$ over $\mathcal{K}$.
In the last  section \ref{sec:theorem} Theorems \ref{thm:1k} and \ref{thm:2} are concluded.

\section{Ordinals for $\Pi^{1}_{N}$-indescribable cardinals}\label{sect:ordinalnotation}
In this section
iterated Skolem hulls $\mathcal{H}_{\alpha,n}(X)$ of sets $X$ of ordinals, ordinals $\Psi_{\kappa,n}\gamma$ for regular ordinals 
$\kappa\,(\mathcal{K}<\kappa\leq I)$ or $\kap=\ome_{1}$,
ordinals $\Psi_{\calk,n}^{\vec{\alp},\Tht}(\alp)$
and classes $Mh_{k,n}(\vec{\alpha})[\Theta]$ are defined for sequences $\vec{\alp}$ of ordinals and finite sets $\Theta$ of ordinals.
It is shown that for each $k<N$ and each $n,m<\omega$,
$\mathcal{K}\in Mh_{k,n}((\omega_{m}(I+1),\ldots,\ome_{m}(I+1))[\emptyset])$ in
${\sf ZF}+(V=L)$ assuming $\mathcal{K}$ is a $\Pi^{1}_{N}$-indescribable cardinal,
where $\ome_{m}(I+1)$ occurs $(N-k)$ times in the class $Mh_{k,n}((\omega_{m}(I+1),\ldots,\ome_{m}(I+1)))[\emptyset]$.

Let $ORD\subset V$ denote the class of ordinals, $ORD^{\varepsilon}\subset V$ and $<^{\varepsilon}$ be $\Delta$-predicates such that
for any transitive and wellfounded model $V$ of $\mbox{{\sf KP}}\omega$,
$<^{\varepsilon}$ is a well ordering of type $\varepsilon_{I+1}$ on $ORD^{\varepsilon}$
for the order type $I$ of the class $ORD$ in $V$.



$<^{\varepsilon}$ is assumed to be a canonical ordering such that 
$\mbox{{\sf KP}}\omega$ proves the fact that $<^{\varepsilon}$ is a linear ordering, and for any formula $\varphi$
and each $n<\omega$,
\begin{equation}\label{eq:trindveps}
\mbox{{\sf KP}}\omega\vdash\forall x(\forall y<^{\varepsilon}x\,\varphi(y)\to\varphi(x)) \to \forall x<^{\varepsilon}\lceil\omega_{n}(I+1)\rceil\varphi(x)
\end{equation}
for the code $\lceil\omega_{n}(I+1)\rceil\in ORD^{\varepsilon}$ of the `ordinal'  $\omega_{n}(I+1)$.

For a definition of $\Delta$-predicates $ORD^{\varepsilon}$ and $<^{\varepsilon}$, 
and a proof of (\ref{eq:trindveps}), cf. \cite{liftupZF}.

In the definition of $ORD^{\varepsilon}$ and $<^{\varepsilon}$,
$I$ with its code $\lceil I\rceil=\langle 1,0\rangle$ 
is {\it intended\/} to denote  the least weakly inaccessible cardinal
above the least $\Pi^{1}_{N}$-indescribable cardinal $\mathcal{K}$,
though we {\it do not assume\/} the existence of weakly inaccessible cardinals above $\mathcal{K}$
anywhere in this paper.
We are working in ${\sf ZF}+(V=L)$ assuming $\mathcal{K}$ is a regular cardinal.

Let
\[
Reg:=\{\ome_{1}\}\cup\{\kap<I: \calk<\kap \mbox{ is regular}\}
\]
while
$Reg^{+}:=Reg\cup\{I\}$.
$\kappa,\lambda,\rho,\pi$ denote elements of $Reg$.
$\kappa^{+}$ denotes the least regular ordinal above $\kappa$.
$\Theta$ denotes finite sets of ordinals$\leq\mathcal{K}$.
$\Theta\subset_{fin} X$ iff $\Theta$ is a finite subset of $X$.
$ORD$ denotes the class of ordinals less than $I$, 
while $ORD^{\varepsilon}$ the class of codes of ordinals less than the next epsilon number $\varepsilon_{I+1}$ to $I$.

For admissible ordinals $\sigma$ and $X\subset L_{\sigma}$, $\mbox{{\rm Hull}}_{\Sigma_{n}}^{\sigma}(X)$ denotes  the $\Sigma_{n}$-Skolem hull of $X$ over $L_{\sigma}$, cf. \cite{liftupZF}.
$F(y)=F^{\Sigma_{n}}(y;\sigma,X)$ denotes the Mostowski collapsing
$F: \mbox{{\rm Hull}}_{\Sigma_{n}}^{\sigma}(X)\leftrightarrow L_{\gamma}$
of $\mbox{{\rm Hull}}_{\Sigma_{n}}^{\sigma}(X)$ for a $\gamma$.
Let $F^{\Sigma_{n}}(\sigma;\sigma,X):=\gamma$.
When $\sigma=I$, we write $F^{\Sigma_{n}}_{X}(y)$ for $F^{\Sigma_{n}}(y;I,X)$.

{\it In what follows up to the last section \ref{sec:theorem},\/} $n\geq 1$ {\it denotes a fixed positive integer\/}.
\\

$Code^{\varepsilon}$ denotes the union of $L_{I}$ and the codes $ORD^{\varepsilon}$ of ordinals$<\varepsilon_{I+1}$.
On $Code^{\varepsilon}$, $\lceil x\rceil \in^{\varepsilon}\lceil y\rceil:\Leftrightarrow x\in y$.
For simplicity let us identify the code $x\in Code^{\varepsilon}$ with
the `set' coded by $x$,
and $\in^{\varepsilon}$ [$<^{\varepsilon}$] is denoted by $\in$ [$<$], resp. when no confusion likely occurs.

Let
$\alpha<\omega_{n+1}(I+1)$ be ordinals, $X\subset L_{I}$ sets, 
$n<\omega$, $0\leq k\leq N$, and
$\kappa\in Reg^{+}$ uncountable regular cardinals$\leq I$.

Define simultaneously 
classes $\mathcal{H}_{\alpha,n}(X)\subset L_{I}\cup\{x\in ORD^{\varepsilon}: x<\omega_{n}(I+1)\}$, 
 and
 ordinals $\Psi_{\kappa,n} \alpha$ as in \cite{liftupZF}.
  We see that $\mathcal{H}_{\alpha,n}(X)$ and $\Psi_{\kappa,n} \alpha$ are 
(first-order) definable as a fixed point in ${\sf ZF}+(V=L)$,
 cf. Proposition \ref{prp:definability}.


\bdf\label{df:Cpsiregularsm}

$\mathcal{H}_{\alpha,n}(X)$ {\rm is the} 
{\rm Skolem hull of} $\{0,\mathcal{K},I\}\cup X$ {\rm under the functions} 
$+,
 \alpha\mapsto\omega^{\alpha}, 
 \kappa\mapsto\kappa^{+}\, (\calk\leq\kappa\in Reg),
(\kap,\gam)\mapsto \Psi_{\kappa,n}\gam\,(\gam<\alp,\kappa\in Reg^{+})$,
{\rm  the Skolem hullings:}
\[
X \mapsto  \mbox{{\rm Hull}}^{I}_{\Sigma_{n}}(X\cap I)
\]
{\rm and the Mostowski collapsing functions}
\[
x=\Psi_{\kappa,n}\gamma\mapsto F^{\Sigma_{1}}_{x\cup\{\kappa\}}\,(\kappa\in Reg)
\]
{\rm and}
\[
x=\Psi_{I,n}\gamma\mapsto F^{\Sigma_{n}}_{x}
\]

\benu
\item
{\rm (Inductive definition of} $\mathcal{H}_{\alpha,n}(X)${\rm ).}
\benu
\item
$\{0,\ome_{1},\mathcal{K},I\}\cup X\subset\mathcal{H}_{\alpha,n}(X)$.

\item
 $x, y \in \mathcal{H}_{\alpha}(X) \Rightarrow x+ y\in \mathcal{H}_{\alpha,n}(X)$,
 {\rm and} $x\in\calh_{\alp,n}(X)\cap\ome_{n}(I+1) \Rarw \omega^{x}\in \mathcal{H}_{\alpha,n}(X)$.
 
\item
$\calk\leq\kappa\in\mathcal{H}_{\alpha,n}(X)\cap Reg \Rightarrow \kappa^{+}\in\mathcal{H}_{\alpha,n}(X)$.

\item
$
\gamma\in \mathcal{H}_{\alpha,n}(X)\cap\alpha
\Rightarrow 
\Psi_{I,n}\gamma\in\mathcal{H}_{\alpha,n}(X)
$.

\item

{\rm If}
$\kappa\in\mathcal{H}_{\alpha,n}(X)\cap Reg$
{\rm and}
$\gamma\in \mathcal{H}_{\alpha,n}(X)\cap\alpha$ 
{\rm then}
$\Psi_{\kappa,n}\gamma\in\mathcal{H}_{\alpha,n}(X)$.




\item

\[
\mbox{{\rm Hull}}^{I}_{\Sigma_{n}}(\mathcal{H}_{\alpha,n}(X)\cap L_{I})\cap Code^{\varepsilon} 
\subset \mathcal{H}_{\alpha,n}(X)
.\]
{\rm Namely for any} $\Sigma_{n}${\rm -formula} $\varphi[x,\vec{y}\,]$ {\rm in the language} $\{\in,=\}$
{\rm and parameters} $\vec{a}\subset \mathcal{H}_{\alpha,n}(X)\cap L_{I}${\rm , if} $b\in L_{I}$,
$(L_{I},\in)\models\varphi[b,\vec{a}\,]$ {\rm and} 
$(L_{I},\in)\models\exists!x\,\varphi[x,\vec{a}\,]$,
{\rm then} $b\in\mathcal{H}_{\alpha,n}(X)$.

\item
{\rm If} $\kappa\in\mathcal{H}_{\alpha,n}(X)\cap Reg$, 
$\gamma\in \mathcal{H}_{\alpha,n}(X)\cap\alpha$, $x=\Psi_{\kappa,n}\gamma\in\mathcal{H}_{\alpha,n}(X)$, 
$\kappa\in\mathcal{H}_{\gamma,n}(\kappa)$
{\rm and}
$\delta\in (\mbox{{\rm Hull}}^{I}_{\Sigma_{1}}(x\cup\{\kappa\})\cup\{I\})\cap\mathcal{H}_{\alpha,n}(X)$,
{\rm then}
$F^{\Sigma_{1}}_{x\cup\{\kappa\}}(\delta)\in\mathcal{H}_{\alpha,n}(X)$.

\item
{\rm If} 
$\gamma\in \mathcal{H}_{\alpha,n}(X)\cap\alpha$, $x=\Psi_{I,n}\gamma\in\mathcal{H}_{\alpha,n}(X)$, 
{\rm and}
$\delta\in (\mbox{{\rm Hull}}^{I}_{\Sigma_{n}}(x)\cup\{I\})\cap\mathcal{H}_{\alpha,n}(X)$,
{\rm then}
$F^{\Sigma_{n}}_{x}(\delta)\in\mathcal{H}_{\alpha,n}(X)$.

\eenu

\item
\[
\calh_{\alp,n}[Y](X):=\calh_{\alp,n}(Y\cup X)
\]
{\rm for sets} $Y\subset L_{I}$.

\item
{\rm (Definition of} $\Psi_{\kappa,n}\alpha${\rm ).}

{\rm Assume} $\kappa\in Reg^{+}${\rm . Then}
\[
\Psi_{\kappa,n}\alpha:=
\min(\{\kappa\}\cup\{\beta<\kappa :  \mathcal{H}_{\alpha,n}[\{\kap\}](\beta)\cap \kappa \subset\beta\})
.\]
\eenu
\edf

Next classes $Mh_{k,n}(\vec{\alp})[\Tht]$ of regular cardinals and
ordinals $\Psi_{\calk,n}^{\vec{\alp},\Tht}(\alp)$ 
are defined.

\bdf\label{df:Lam}

\benu
\item
{\rm Let} $\vec{\alp}=(\alp_{0},\ldots,\alp_{m-1})$ {\rm be a sequence of ordinals.}

 \benu
 \item
 length $lh(\vec{\alp}):=m$, components $\vec{\alp}(i):=\alp_{i}$ {\rm and}
 end segments $\vec{\alp}[i]:=(\alp_{i},\ldots,\alp_{m-1})$ {\rm for} $i<lh(\vec{\alp})$.

\item
{\rm The set of components}
\[
K(\vec{\alp}):=\{\vec{\alp}(i): i<lh(\vec{\alp}) \}=\{\alp_{0},\ldots,\alp_{m-1}\}
.\]

 \item
{\rm Sequences consisting of a single element} $(\alp)$ {\rm is identified with the ordinal} $\alp$,
{\rm and} $\emptyset$ {\rm denotes the} empty sequence.

 \item

{\rm For sequence of ordinals} $\vec{\nu}$ {\rm of the same length}, $lh(\vec{\nu})= lh(\vec{\alp})$,
\[
\vec{\nu}<\vec{\alp}   :\Lrarw  \fal i<lh(\vec{\alp})[\vec{\nu}(i)<\vec{\alp}(i)]
\]

 \item
 {\rm For ordinals} $\bet$,
\[
\vec{\alp}\leq\bet :\Lrarw \fal\alp\in K(\vec{\alp})(\alp\leq\bet)
.\]

 \item
{\rm For sequence of ordinals} $\vec{\nu}=(\nu_{0},\ldots,\nu_{m-1})$ {\rm of the same length}, $lh(\vec{\nu})=lh(\vec{\alp})$
{\rm and} $i<lh(\vec{\alp})$
\[
(\vec{\nu}\bullet\vec{\alp})[i]  := (\vec{\nu}(i))* \vec{\alp}[i+1] =(\nu_{i},\alp_{i+1},\ldots,\alp_{m-1})
.\]
{\rm Note that} $lh((\vec{\nu}\bullet\vec{\alp})[i])=lh(\vec{\alp})-i$.
 \eenu

\item
{\rm For} $A\subset ORD$, {\rm limit ordinals} $\alp$ {\rm and} $i\geq 0$
\[
\alp\in M_{i}(A) :\Lrarw  A \mbox{ {\rm is} } \Pi^{1}_{i}\mbox{{\rm -indescribable in} } \alp
(\Lrarw A \mbox{ {\rm is} } (i+1)\mbox{{\rm -stationary in }} \alp).
\]

\eenu
\edf

\bdf\label{df:psivecalp}

\benu

\item

{\rm For sequence of ordinals} $\vec{\nu},\vec{\alp}$ {\rm in the same length,  let}
\[
\vec{\nu}\in \mathcal{H}_{\vec{\nu},n}[\Theta](\pi)\cap\vec{\alp}
 :\Lrarw 
\fal i<lh(\vec{\nu})[\vec{\nu}(i)\in\calh_{\vec{\nu}(i),n}[\Tht](\pi)\cap\vec{\alp}(i)]
\]


\item {\rm (Definition of $Mh_{k,n}(\vec{\alp})[\Tht]$)}
{\rm First let for the empty sequence} $\emptyset$
\[
\calk\in Mh_{N,n}(\emptyset)[\Tht] :\Lrarw 
\calk\in M_{N} 
\Lrarw \calk \mbox{ {\rm is }} \Pi^{1}_{N}\mbox{{\rm -indescribable}}.
\]

{\rm The classes} $Mh_{k,n}(\vec{\alp})[\Tht]$ {\rm are defined for} $n<\omega$, $0\leq k< N$,
 $\Tht\subset_{fin}(\calk+1)$, {\rm and sequences of ordinals} $\vec{\alp}$ {\rm such that}
$lh(\vec{\alp})=N-k$.
\beqnarr
&& \pi\in Mh_{k,n}(\vec{\alp})[\Tht]  :\Lrarw 
\pi \mbox{ {\rm is a regular cardinal}}\leq\calk \mbox{ {\rm such that}}
\label{eq:dfMhkh}
\\
&&
\fal \vec{\nu}\in\calh_{\vec{\nu},n}[\Tht\cup\{\pi\}](\pi)\cap\vec{\alp}
[
\pi\in M_{k}(
\bigcap_{i<lh(\vec{\alp})}
Mh_{k+i,n}((\vec{\nu}\bullet\vec{\alp})[i] )[\Tht\cup\{\pi\}]
)
]
\nonumber
\eeqnarr

 \item
{\rm We say that the class} 

\beqn
\label{eq:Hdf}
\{\rho\in \bigcap_{i<lh(\vec{\alp})}
Mh_{k+i,n}((\vec{\nu}\bullet\vec{\alp})[i] )[\Tht\cup\{\pi\}]\cap\pi :
\calh_{\gam,n}[\Tht\cup\{\pi\}](\rho)\cap\pi\subset\rho \}
\eeqn
{\rm is the} resolvent class {\rm for} $\pi\in Mh_{k,n}(\vec{\alp})[\Tht]$
{\rm with respect to} $\vec{\nu}$ {\rm and} $\gam$.

\item {\rm (Definition of $\Psi_{\calk,n}^{\vec{\alp},\Tht}(\alp)$)}

{\rm For} $\vec{\alp}\leq\alp$ {\rm with} $lh(\vec{\alp})=N$, {\rm let}
{\small
\beqn\label{eq:Psivec}
\Psi_{\calk,n}^{\vec{\alp},\Tht}(\alp):=\min(\{\calk\}\cup
\{\pi\in \bigcap_{i<N}Mh_{i,n}(\vec{\alp}[i])[\Tht\cup\{\pi\}]\cap\calk: 
  \calh_{\alp,n}[\Tht](\pi)\cap\calk\subset\pi 
\})
\eeqn
}

\item
\[
Mh_{k,n}(\vec{\alp}):=Mh_{k,n}(\vec{\alp})[\emptyset]
.\]

\eenu

\edf

$x,y,z,\ldots$ range over sets in $L_{I}$,
$\alp,\bet,\gam,\ldots$ range over $ORD^{\veps}$,
$\vec{\alp},\vec{\nu},\ldots$ range over finite sequences over $ORD^{\veps}$.
$\vphi,\tau$ denote formulae.

The following Proposition \ref{prp:definability} is easy to see.

\begin{proposition}\label{prp:definability}
Each of 
$x=\mathcal{H}_{\alpha,n}(\beta)\, (\alpha\in ORD^{\varepsilon},\beta<^{\varepsilon}I)$,
$\beta=\Psi_{\kappa,n}\alpha\,(\kappa\in R^{+})$,
$x\in Mh_{k,n}(\vec{\alp})[\Theta]$ and $\bet=\Psi^{\vec{\alp},\Tht}_{\calk,n}(\alp)$
is a $\Sigma_{n+1}$-predicate as fixed points in ${\sf ZF}+(V=L)$. 
\end{proposition}
\bprf 
Let us examine the definability of $x\in Mh_{k,n}(\vec{\alp})[\Theta]$.

Let $\alp=\max K(\vec{\alp})$, and $m=lh(\vec{\alp})$.
Then $\pi\in Mh_{k,n}(\vec{\alp})[\Tht]$ iff $\pi\leq\calk$ is regular and
there exist a set
$y=\calh_{\alp,n}[\Tht\cup\{\pi\}](\pi)$ and a function $\{x_{\nu}\}_{\nu\in y}$ such that
$\fal\nu\in y[x_{\nu}=\calh_{\nu,n}[\Tht\cup\{\pi\}](\pi)]$ and for any sequence $\vec{\nu}\in{}^{m}y$(, i.e., $\fal i<m=lh(\vec{\nu})(\vec{\nu}(i)\in y)$)
with $\vec{\nu}<\vec{\alp}$,
if $\fal i<lh(\vec{\nu})(\vec{\nu}(i)\in x_{\vec{\nu}(i)})$, then
for any $\Pi^{1}_{k}(\pi)$-sentence $\tht$ true on $L_{\pi}$,
there exists a $\sig<\pi$ such that $\tht$ holds in $L_{\sig}$ and
$ \sig\in \bigcap_{i<lh(\vec{\alp})}Mh_{k+i,n}((\vec{\nu}\bullet\vec{\alp})[i])[\Tht\cup\{\pi\}]$.
\eprf

\begin{proposition}\label{prp:clshull}

\benu

\item\label{prp:Mh3}
Suppose $lh(\vec{\bet})=lh(\vec{\alp})$ and $\fal i<lh(\vec{\alp})(\vec{\bet}(i)\leq\vec{\alp}(i))$.
Then $\pi\in Mh_{k,n}(\vec{\alpha})[\Theta]  \Rightarrow \pi\in Mh_{k,n}(\vec{\bet})[\Theta]$.

\item\label{prp:Mh2}
$\calk\in M_{N-1}(Mh_{N-1,n}(\alpha)[\Theta\cup\{\calk\}])\Rightarrow \calk\in Mh_{N-1,n}(\alpha)[\Theta\cup\{\calk\}]$.

\item\label{prp:L4.10.1.1a}
$
Mh_{k,n}(\vec{\alp})[\Tht\cup\{\kap\}]\subset Mh_{k,n}(\vec{\alp})[\Tht]
$.

\item\label{prp:L4.10.1.1aa}
$
\rho\in Mh_{k,n}(\vec{\alp})[\Tht]\Lrarw \rho\in Mh_{k,n}(\vec{\alp})[\Tht\cup\{\rho\}]
$.

\item\label{prp:L4.10.1.1c}
$\kap\in Mh_{k,n}(\vec{\alp})[\Tht] \Lrarw \kap\in Mh_{k,n}(\vec{\alp})[\Tht\cup\{\bet\}]$
for any $\bet<\kap$.

\item\label{prp:L4.10.1.1d}
$\Tht_{1}\subset_{fin}\pi \spand \pi\in M_{k}(Mh_{k,n}(\vec{\alp})[\Tht_{0}])\Rarw
\pi\in M_{k}(Mh_{k,n}(\vec{\alp})[\Tht_{0}\cup\Tht_{1}])$.

\eenu
\end{proposition}
\bprf
\ref{prp:clshull}.\ref{prp:Mh3}. This is seen from Definition \ref{df:psivecalp}, (\ref{eq:dfMhkh}).
\\

\noindent
\ref{prp:clshull}.\ref{prp:Mh2}. 
Suppose $\calk\in M_{N-1}(Mh_{N-1,n}(\alpha)[\Theta\cup\{\calk\}])$ and $\bet<\alp$.
Then by Proposition \ref{prp:clshull}.\ref{prp:Mh3} we have
$\calk\in M_{N-1}(Mh_{N-1,n}(\bet)[\Theta\cup\{\calk\}])$, and hence
$\calk\in Mh_{N-1,n}(\alpha)[\Theta\cup\{\calk\}]$.
\\

\noindent
\ref{prp:clshull}.\ref{prp:L4.10.1.1a}.
$\lam\in Mh_{k,n}(\vec{\alp})[\Tht\cup\{\kap\}]\Rarw \lam\in Mh_{k,n}(\vec{\alp})[\Tht]$ is seen easily by induction on $\lam$ 
from 
the monotonicity of the operators $M_{k}$,
$X\subset Y \Rarw M_{k}(X)\subset M_{k}(Y)$, and
$\calh_{\nu,n}[\Tht\cup\{\lam\}]\subset\calh_{\nu,n}[\Tht\cup\{\kap,\lam\}]$.
\\

\noindent
\ref{prp:clshull}.\ref{prp:L4.10.1.1aa}. This is seen from Definition \ref{df:psivecalp}, (\ref{eq:dfMhkh})
and Proposition \ref{prp:clshull}.\ref{prp:L4.10.1.1a}.
\\

\noindent
\ref{prp:clshull}.\ref{prp:L4.10.1.1c}.
By induction on $\kap$.
Let $\bet<\kap$.
We have
$\vec{\nu}\in\calh_{\vec{\nu},n}[\Tht\cup\{\kap\}](\kap)$ iff
$\vec{\nu}\in\calh_{\vec{\nu},n}[\Tht\cup\{\kap,\bet\}](\kap)$
since $\bet<\kap$.
It suffices to show for such a $\vec{\nu}<\vec{\alp}$,
$\kap\in M_{k}(\bigcap_{i<lh(\vec{\alp})}Mh_{k+i,n}((\vec{\nu}\bullet\vec{\alp})[i])[\Tht\cup\{\kap\}])$ iff
$\kap\in M_{k}(\bigcap_{i<lh(\vec{\alp})}Mh_{k+i,n}((\vec{\nu}\bullet\vec{\alp})[i])[\Tht\cup\{\kap,\bet\}])$.
By IH for any $\rho$ with $\bet<\rho<\kap$,
$\rho\in Mh_{k+i,n}((\vec{\nu}\bullet\vec{\alp})[i])[\Tht\cup\{\kap\}]$ iff
$\rho\in Mh_{k+i,n}((\vec{\nu}\bullet\vec{\alp})[i])[\Tht\cup\{\kap,\bet\}]$.
This shows the lemma since for any $\bet<\kap$ and any $X\subset ORD$,
$\kap\in M_{k}(X)$ iff $\kap\in M_{k}(X\setm(\bet+1))$.
\\

\noindent
\ref{prp:clshull}.\ref{prp:L4.10.1.1d}.
This is seen from Proposition \ref{prp:clshull}.\ref{prp:L4.10.1.1c}.

\eprf

Let $A_{n}(\alpha)$ denote the conjunction of 
$\forall\beta<^{\varepsilon}I\exists ! x[x=\mathcal{H}_{\alpha,n}(\beta)]$,
$\forall\kappa\in R^{+}\forall x[x=\mathcal{H}_{\alpha,n}(\kappa) \to \exists!\beta<\kap(\beta=\Psi_{\kappa,n}\alpha)]$ and 
$\fal\vec{\alp}\fal k< N\forall\Theta\subset_{fin}(\mathcal{K}+1)(\max K(\vec{\alp})\leq\alp \spand 
lh(\vec{\alp})+k=N \to
\exists ! x (x=Mh_{k,n}(\vec{\alp})[\Theta]) \spand
\exi!\bet(\bet=\Psi^{\vec{\alp},\Tht}_{\calk,n}(\alp)))$.
\\

$card(x)$ denotes the cardinality of sets $x$.

\begin{lemma}\label{lem:welldefinedness}
For each $n,m<\omega$, ${\sf ZF}+(V=L)$ proves the followings.
\benu
\item\label{lem:welldefinedness.0}
$y=\mathcal{H}_{\alpha,n}(x) \to card(y)\leq\max\{card(x),\aleph_{0}\}$.

\item\label{lem:welldefinedness.2}
$\forall\alpha<^{\varepsilon}\omega_{m}(I+1)\, A_{n}(\alpha)
$.

\item\label{lem:welldefinedness.1}
If $\mathcal{K}$ is $\Pi^{1}_{N}$-indescribable and $\Tht\subset_{fin}(\calk+1)$, then
$\mathcal{K}\in Mh_{N-1,n}(\omega_{m}(I+1))[\Theta]\cap M_{N-1}(Mh_{N-1,n}(\omega_{m}(I+1))[\Theta])
$ for each $m<\ome$.

\item\label{lem:welldefinedness.3}
$\pi\in Mh_{k,n}(\vec{\alp})[\Tht]$ is a $\Pi^{1}_{k+1}$-class on $L_{\pi}$ uniformly for weakly inaccessible cardinals $\pi\leq\calk$.
This means that for each $k,n$ 
 there exist a $\Pi^{1}_{k+1}$-formula $m_{k,n}(R_{0},R_{1})$ and
some $r_{i}\subset L_{\calk}$ depending on $\vec{\alp},\Tht$ such that
$\pi\in Mh_{k,n}(\vec{\alp})[\Tht]$ iff $L_{\pi}\models m_{k,n}(r_{0}\cap L_{\pi},r_{1}\cap L_{\pi})$ 
for any weakly inaccessible cardinals $\pi\leq\calk$.
\eenu
\end{lemma}
\bprf

\ref{lem:welldefinedness}.\ref{lem:welldefinedness.2}.
We show that $A_{n}(\alpha)$ is progressive, i.e.,
$\forall\alpha<^{\varepsilon}\omega_{m}(I+1)[\forall\gamma<^{\varepsilon}\alpha\, A_{n}(\gamma) \to A_{n}(\alpha)]$.

Assume $\forall\gamma<^{\varepsilon}\alpha\, A_{n}(\gamma)$ and $\alpha<^{\varepsilon}\omega_{m}(I+1)$.
$\forall\beta<^{\varepsilon}I\exists ! x[x=\mathcal{H}_{\alpha,n}(\beta)]$
 follows from IH and the Replacement.

Next assume $\kappa\in Reg^{+}$.
Then $\exists!\beta<\kap(\beta=\Psi_{\kappa,n}\alpha)$ follows from the regularity of $\kappa$.

$\max K(\vec{\alp})\leq\alp \to \exists ! x[x=Mh_{k,n}(\vec{\alp})[\Theta]]$ is easily seen from IH,
Proposition \ref{prp:definability} and Separation\, $Mh_{k,n}(\vec{\alp})[\Theta]\subset(\calk+1)$.
$\exi!\bet(\bet=\Psi^{\vec{\alp},\Tht}_{\calk,n}(\alp))$ follows from this.
\\

\noindent
\ref{lem:welldefinedness}.\ref{lem:welldefinedness.1}.
Suppose $\mathcal{K}$ is $\Pi^{1}_{N}$-indescribable.
We show the following $B_{n}(\alpha)$ is progressive in $\alpha\in ORD^{\veps}$:
\beqnarrs
B_{n}(\alp) & :\Lrarw &
\fal\Tht\subset_{fin}(\calk+1)
[\alpha\in \mathcal{H}_{\alpha,n}[\Theta\cup\{\calk\}](\mathcal{K})
\\
&&
\to \mathcal{K}\in Mh_{N-1,n}(\alpha)[\Theta] \cap M_{N-1}(Mh_{N-1,n}(\alpha)[\Theta] )]
.\eeqnarrs

Suppose $\forall\xi<^{\varepsilon}\alpha\, B_{n}(\xi)$, $\Tht\subset_{fin}(\calk+1)$, 
 and $\alpha\in \mathcal{H}_{\alpha,n}[\Theta\cup\{\calk\}](\mathcal{K})$.
We have to show that $Mh_{N-1,n}(\alpha)[\Theta]$ is $\Pi^{1}_{N-1}$-indescribable in $\mathcal{K}$.
$\mathcal{K}\in Mh_{N-1,n}(\alpha)[\Theta]$ follows from 
$\mathcal{K}\in M_{N-1}(Mh_{N-1,n}(\alpha)[\Theta])$, cf. Proposition \ref{prp:clshull}.\ref{prp:Mh2}.

Let $\sig(X,m)$ be a universal $\Pi^{1}_{N-1}$-formula,
and assume that $L_{\calk}\models\sig(C_{0},m_{0})$ for a subset $C_{0}$ of $\calk$ and an $m_{0}<\ome$.


Since $\forall\pi\leq\mathcal{K}[card(\mathcal{H}_{\alpha,n}[\Theta\cup\{\pi\}](\pi))\leq\pi]$, 
pick an injection $f:\mathcal{H}_{\alpha,n}[\Theta\cup\{\mathcal{K}\}](\mathcal{K})\to \mathcal{K}$ so that
$f"\mathcal{H}_{\alpha,n}[\Theta\cup\{\pi\}](\pi)\subset\pi$ for any weakly inaccessibles $\pi\leq\mathcal{K}$.

Let $R_{0}=\{f(\alpha)\}$, $R_{1}=C_{0}$, 
$R_{2}=\{f(\xi) :\xi\in \mathcal{H}_{\xi,n}[\Theta\cup\{\calk\}](\mathcal{K})\cap\alpha\}$,
$R_{3}=\bigcup\{(Mh_{N-1,n}(\xi)[\Theta\cup\{\pi\}]\cap\mathcal{K})\times\{f(\pi)\}\times\{f(\xi)\} :\xi\in \mathcal{H}_{\xi,n}[\Theta\cup\{\calk\}](\mathcal{K})\cap\alpha, \pi\leq\mathcal{K}\}$,
and $R_{4}=\{(f(\beta),f(\gamma)): \{\beta,\gamma\}\subset\mathcal{H}_{\alpha,n}[\Theta\cup\{\mathcal{K}\}](\mathcal{K}), \beta<\gamma\}$.


By IH we have $\forall\xi\in \mathcal{H}_{\xi,n}[\Theta\cup\{\calk\}](\mathcal{K})\cap\alpha[\mathcal{K}\in M_{N-1}(Mh_{N-1,n}(\xi)[\Theta\cup\{\calk\}])]$.
Hence $\langle L_{\mathcal{K}},\in, R_{i}\rangle_{i\leq 4}$
enjoys a $\Pi^{1}_{N}$-sentence saying that
$\mathcal{K}$ is weakly inaccessible, $R_{0}\neq\emptyset$, 
$\sig(R_{1},m_{0})$ and
{\small
\[
\varphi:\Leftrightarrow 
\fal m\in\ome\forall C\subset ORD \forall x, y \exi a[
R_{2}(x) \land \theta(R_{4},y) \land \sig(C,m)
\to R_{3}(a,y,x) \land L_{a}\models\sig(C\cap a,m)]
\]
}
where 
$\theta(R_{4},y)$ is a $\Sigma^{1}_{1}$-formula such that for any $\pi\leq\mathcal{K}$
\[
L_{\pi}\models \theta(R_{4},y) \Leftrightarrow y=f(\pi)
.\]

By the $\Pi^{1}_{N}$-indescribability of $\mathcal{K}$, pick a $\pi<\mathcal{K}$ such that
$\langle L_{\pi},\in, R_{i}\cap L_{\pi}\rangle_{i\leq 4}$ enjoys the $\Pi^{1}_{N}$-sentence.

We claim that $\pi\in Mh_{N-1,n}(\alpha)[\Theta]$ and $L_{\pi}\models\sig(C_{0}\cap\pi,m_{0})$.

$\pi$ is weakly inaccessible, $f(\alpha)\in \pi$ and 
$\sig(C_{0}\cap\pi,m_{0})$ hold in $L_{\pi}$.
We can assume that
$\calh_{\alp,n}[\{\calk\}](\pi)\cap\calk\subset\pi$
since $\{\pi<\calk:\calh_{\alp,n}[\{\calk\}](\pi)\cap\calk\subset\pi\}$ is club in $\calk$.

It remains to see 
$\forall\xi\in \mathcal{H}_{\xi,n}[\Theta\cup\{\pi\}](\pi)\cap\alpha[\pi\in M_{N-1}(Mh_{N-1,n}(\xi)[\Theta\cup\{\pi\})]$.
This follows from the fact that
$\varphi$ holds in $\langle L_{\pi},\in, R_{i}\cap L_{\pi}\rangle_{i\leq 4}$, and
$\forall\xi\in \mathcal{H}_{\xi,n}[\Theta\cup\{\pi\}](\pi)\cap\alpha(f(\xi)\in \pi)$
by $f" \mathcal{H}_{\alpha,n}[\Theta\cup\{\pi\}](\pi)\subset\pi$ and 
$\mathcal{H}_{\xi,n}[\Theta\cup\{\pi\}](\pi)\subset\mathcal{H}_{\xi,n}[\Theta\cup\{\calk\}](\calk)$.

Thus $\mathcal{K}\in M_{N-1}(Mh_{N-1,n}(\alpha)[\Theta])$.
\\

\noindent
\ref{lem:welldefinedness}.\ref{lem:welldefinedness.3}.
This is seen as in Lemma \ref{lem:welldefinedness}.\ref{lem:welldefinedness.1} using a universal $\Pi^{1}_{k}$-formula
and an injection $f:\mathcal{H}_{\alpha,n}[\Theta\cup\{\mathcal{K}\}](\mathcal{K})\to \mathcal{K}$ 
for $\alp=\max K(\vec{\alp})$ so that
$f"\mathcal{H}_{\alpha,n}[\Theta\cup\{\pi\}](\pi)\subset\pi$ for any weakly inaccessibles $\pi\leq\mathcal{K}$.
\eprf

\blem\label{lem:stepdown}
For each $m<\ome$, ${\sf ZF}+(V=L)$ proves the following.

Let $\pi\in Mh_{k+1,n}(\vec{\alp})[\Tht]$.
Then for any ordinals $\gam<\ome_{m}(I+1)$, 
$\pi\in Mh_{k,n}((\gam)*\vec{\alp})[\Tht]$.


\elem
\bprf
Suppose $\pi\in Mh_{k+1,n}(\vec{\alp})[\Tht]$.
We show $
\pi\in Mh_{k,n}((\gam)*\vec{\alp})[\Tht]$
by induction on $\gam$.

Let $X=\bigcap_{i<lh(\vec{\alp})} Mh_{k+1+i,n}((\vec{\nu}\bullet\vec{\alp})[i])[\Tht\cup\{\pi\}]$ and
$Y=Mh_{k,n}((\del)*\vec{\alp})[\Tht\cup\{\pi\}]$
for a sequence $\vec{\nu}\in\calh_{\vec{\nu},n}[\Tht\cup\{\pi\}](\pi)\cap \vec{\alp}$ and
an ordinal $\del\in\calh_{\del,n}[\Tht\cup\{\pi\}](\pi)\cap\gam$.
We show
$\pi\in M_{k}(X\cap Y)$, cf. the definition (\ref{eq:dfMhkh}).

The assumption $\pi\in Mh_{k+1,n}(\vec{\alp})[\Tht]$
yields $\pi\in M_{k+1}(X)$.
On the other side IH yields 
$\pi\in Mh_{k,n}((\del)*\vec{\alp})[\Tht]$.
By Proposition \ref{prp:clshull}.\ref{prp:L4.10.1.1aa} we have $\pi\in Y$.

Since $\pi\in Y$ is a $\Pi^{1}_{k+1}$-sentence on $L_{\pi}$ by Lemma \ref{lem:welldefinedness}.\ref{lem:welldefinedness.3},
$\pi\in M_{k+1}(X)$
yields
$\pi\in M_{k+1}(X\cap Y)$, a fortiori 
$\pi\in M_{k}(X\cap Y)$.
Thus 
$\pi\in Mh_{k,n}((\gam)*\vec{\alp})[\Tht]$ is shown.


\eprf

\bth\label{th:KM1}

For each $m<\ome$ and $k<N$, ${\sf ZF}+(V=L)$ proves the following.

Let $\mathcal{K}$ be a $\Pi^{1}_{N}$-indescribable cardinal.

Then 
$\mathcal{K}\in Mh_{k,n}(\vec{\bet})[\emptyset]$,
$\mathcal{K}\in Mh_{k,n}(\vec{\bet})[\Tht]$ and $\exi\kap<\calk(\kap=\Psi_{\calk,n}^{\vec{\alp},\Tht}(\alp))$
for any finite set $\Tht\subset_{fin}(\calk+1)$ and any $\vec{\bet}, \alp,\vec{\alp}$ with $lh(\vec{\bet})=N-k$,
$K(\vec{\bet})<\ome_{m}(I+1)$ and
$\vec{\alp}\leq\alp<\ome_{m}(I+1)$.
\end{theorem}
\bprf
By Lemma  \ref{lem:welldefinedness}.\ref{lem:welldefinedness.3} 
$\bigcap_{k<N}Mh_{k,n}(\vec{\alp}[k])[\Tht])$ is a $\Pi^{1}_{N}$-class.
Hence if a $\Pi^{1}_{N}$-indescribable $\calk\in\bigcap_{k<N}Mh_{k,n}(\vec{\alp}[k])[\Tht])$, 
then $\calk\in M_{N}(\bigcap_{k<N}Mh_{k,n}(\vec{\alp}[k])[\Tht]))$,
a fortiori, $\calk\in M_{0}(\bigcap_{k<N}Mh_{k,n}(\vec{\alp}[k])[\Tht]))$.
Since $\{\kap<\calk:\calh_{\alp,n}[\Tht](\kap)\cap\calk\subset\kap\}$ is club in $\calk$,
$\mathcal{K}\in M_{0}(\bigcap_{k<N}Mh_{k,n}(\vec{\alp}[k])[\Tht]))$ 
yields $\exi\kap<\calk(\kap=\Psi_{\calk,n}^{\vec{\alp},\Tht}(\alp))$ for $\vec{\alp}\leq\alp$.

By metainduction on $N-k$ using Lemmata \ref{lem:welldefinedness}.\ref{lem:welldefinedness.1}
and \ref{lem:stepdown}
we see
$\mathcal{K}\in Mh_{k,n}(\vec{\bet})[\emptyset]$ or equivalently 
$\mathcal{K}\in Mh_{k,n}(\vec{\bet})[\{\calk\}]$.
Proposition \ref{prp:clshull}.\ref{prp:L4.10.1.1d} yields 
$\mathcal{K}\in Mh_{k,n}(\vec{\bet})[\Tht\cup\{\calk\}]$ for any $\Tht<\calk$.
\eprf

\section{A theory for $\Pi^{1}_{N}$-indescribable cardinal}\label{sect:Ztheory}

In this section the set theory ${\sf ZF}+(V=L)+(\mathcal{K}\mbox{ is } \Pi^{1}_{N}\mbox{-indescribable)}$
is paraphrased to another set theory $\mbox{{\rm T}}_{N}(\mathcal{K},I)$.
Since our formulation of  $\mbox{{\rm T}}_{N}(\mathcal{K},I)$ is the same as in \cite{liftupK, liftupZF},
let me define it briefly.

Let $\mathcal{K}$ be the least $\Pi^{1}_{N}$-indescribable cardinal for a fixed positive integer $N\geq 1$, 
and $I>\mathcal{K}$ (be intended to denote) 
the least weakly inaccessible cardinal above $\mathcal{K}$.
Again note that we do not assume that such an $I$ exists.
$\kappa,\lambda,\rho$ ranges over uncountable regular ordinals $\sig$ such that $\mathcal{K}<\sig< I$
or $\sig=\ome_{1}$.

In the following Definition \ref{df:regext}, 
$\mbox{{\rm Hull}}^{I}_{\Sigma_{n}}(x)$ denotes the $\Sig_{n}$-Skolem hull of the set $x$ on the universe $L_{I}=L$,
and $\mbox{{\rm Hull}}^{I}_{\Sigma_{n}}(x)\ni z\mapsto F^{\Sigma_{n}}_{x}(z)\in L_{y}$ its Mostowski collapsing map
with $y=F^{\Sigma_{n}}_{x}(I)$.
Then
the predicate $P$ is intended to denote the relation
\[
P(\lambda,x,y)  \Leftrightarrow  x=F^{\Sigma_{1}}_{x\cup\{\lambda\}}(\lambda)
\,\&\,
y=F^{\Sigma_{1}}_{x\cup\{\lambda\}}(I):=rng(F^{\Sigma_{1}}_{x\cup\{\lambda\}})\cap ORD
\]
and the predicate $P_{I,n}(x)$ is intended to denote the relation
\[
P_{I,n}(x)  \Leftrightarrow  x=F^{\Sigma_{n}}_{x}(I)
.\]

\bdf\label{df:pi1kfml}
\benu
\item\label{df:pi1kfml.1}
{\rm For a formula} $\varphi$ {\rm and a set} $x$,
$\varphi^{x}$ {\rm denotes the result of restricting every unbounded quantifier}
$\exists z,\forall z$ {\rm in} $\varphi$ {\rm to} $\exists z\in x, \forall z\in x$.

\item\label{df:pi1kfml.2}
{\rm For natural numbers} $k\geq 0$, {\rm a} $\Pi^{1}_{k}$-formula
{\rm is a formula obtained from a (first-order) formula} $\varphi[\vec{X}\,]$ {\rm  in the language} $\{\in\}\cup\vec{X}$ 
{\rm with unary predicates} $\vec{X}$
{\rm by applying alternating second-order quantifiers at most} $k${\rm -times},
$\fal X_{k}\exi X_{k-1}\cdots Q X_{1}\,\vphi[\vec{X}]$.

$\Sig^{1}_{k}$-formulae {\rm are defined dually.}

\item\label{df:pi1kfml.3}
{\rm Let} $\{X_{i}\}_{i<\ome}$ {\rm be the list of second-order variables, and} $\{x_{i}\}_{i<\ome}$
{\rm the list of first-order variables.}
{\rm Also let} $\vphi(X_{0})$ {\rm be a} $\Pi^{1}_{k}${\rm -formula possibly with a unary predicate} $X_{0}$
{\rm such that for any second-order} $X_{i}$ {\rm occurring in} $\vphi(X_{0})$, {\rm the first-order} $x_{i}$
{\rm does not occur in it}.

{\rm For ordinals} $\alp$ {\rm and natural numbers} $k\geq 0$, {\rm a} $\Pi^{1}_{k}(\alp)$-formula
{\rm (in the language $\{\in\}$, cf. Definition \ref{df:fmlclasses} for the class $\Pi^{1}_{k}(\alp)$ in an expanded language)}
{\rm is a formula obtained from such a formula} $\vphi(X_{0})$ 
{\rm by replacing} $X_{0}(t)$ {\rm by} $t\in x_{0} \land t\in \alp$, {\rm replacing each second-order variable}
$X_{i}(t)$ {\rm occurring in the matrix of} $\vphi(X_{0})$ {\rm by} $t\in x_{i}\land t\in\alp$,
{\rm restricting every second-order quantifier} $\exi X_{i}, \fal X_{i}$ {\rm in} $\vphi$ {\rm to}
$\exi x_{i}\subset \alp, \fal x_{i}\subset \alp$, {\rm and restricting every unbounded quantifier}
$\exists z,\forall z$ {\rm in} $\varphi$ {\rm to} $\exists z\in \alp, \forall z\in \alp$.

\item\label{df:pi1kfml.4}
{\rm For ordinals} $\alp,\bet$ {\rm and a} $\Pi^{1}_{k}(\alp)${\rm -formula} $\vphi\equiv\vphi(y\cap \alp)$,
$\vphi^{(\bet,\alp)}$ {\rm denotes a} $\Pi^{1}_{k}(\bet)${\rm -formula obtained from} $\vphi$
{\rm by replacing} $z\in y\land z\in \alp$ {\rm by} $z\in y \land z\in \bet$,
{\rm restricting every second-order quantifier} $\exi x\subset \alp, \fal x\subset \alp$ {\rm in} $\vphi$ {\rm to}
$\exi x\subset \bet, \fal x\subset \bet$, {\rm and restricting every bounded quantifier}
$\exists z\in \alp,\forall z\in \alp$ {\rm in} $\varphi$ {\rm to} $\exists z\in \bet, \forall z\in \bet$.

\eenu
\end{definition}

\begin{definition}\label{df:regext}
$\mbox{{\rm T}}_{N}(\mathcal{K},I,n)$ {\rm denotes the set theory defined as follows.}
\benu
\item
{\rm Its language is} $\{\in, P,P_{I,n},Reg,\mathcal{K},\ome_{1}\}$ {\rm for a ternary predicate} $P${\rm  , unary predicates} $P_{I,n}$
{\rm and} $Reg${\rm , and individual constants} $\mathcal{K}$ {\rm and} $\ome_{1}$.

\item
{\rm Its axioms are obtained from those of Kripke-Platek set theory with the axiom of infinity} $\mbox{{\rm KP}}\omega$
{\rm in the expanded language,
the axiom of constructibility,}
$V=L$
{\rm together with the following axiom schemata:}
\benu
 \item
 $Reg(\ome_{1})$,
$(Reg(\kappa) \to \kappa\in ORD\land(\kap=\ome_{1}\lor\kap>\calk))$, 
$(Reg(\kappa) \to a\in ORD\cap\kappa \to \exists x, y\in ORD\cap\kappa[a<x\land P(\kappa,x,y)])$, {\rm and}
 $(P(\kappa,x,y) \to \{x,y\}\subset ORD \land  x<y<\kappa \land Reg(\kappa))$ {\rm and}  
$(P(\kappa,x,y) \to a\in ORD\cap x \to \varphi[\kappa,a] \to \varphi^{y}[x,a])$
{\rm for any} $\Sigma_{1}${\rm -formula} $\varphi$ {\rm in the language} $\{\in\}$.

$(\forall x\in ORD\, \exists y[x\geq\mathcal{K} \to y>x\land Reg(y)])$.

 $(P_{I,n}(x)\to x\in ORD)$ {\rm and} $(P_{I,n}(x) \to a\in ORD\cap x \to \varphi[a] \to \varphi^{x}[a])$
{\rm for any} $\Sigma_{n}${\rm -formula} $\varphi$ {\rm in the language} $\{\in\}$,
{\rm and} $(\mathcal{K}<a\in ORD \to \exists x\in ORD[a<x\land P_{I,n}(x)])$.

 \item
 $(\mathcal{K}>\ome_{1})$ 
  {\rm and}

\begin{equation}\label{eq:Kord}
\fal x\subset\calk[\vphi(x)\to\exi \rho<\calk(
\vphi^{(\rho,\calk)}(x\cap\rho))]
\end{equation}
{\rm where} $\vphi(x)$ {\rm is a} $\Pi^{1}_{N}(\calk)${\rm -formula in the language $\{\in\}$.}

{\rm Note that}  `$\mathcal{K}$ {\rm is regular}' , {\rm i.e.,} 
$\forall\alpha<\mathcal{K}\forall f\in{}^{\alpha}\mathcal{K}\exists\beta<\mathcal{K}(f"\alpha\subset\beta)$ {\rm follows from (\ref{eq:Kord}).}
 \eenu

\eenu
\end{definition}

The following Lemma \ref{lem:regularset} is seen as in \cite{liftupK, liftupZF}.

\begin{lemma}\label{lem:regularset}
$\mbox{{\rm T}}_{N}(\mathcal{K},I):=\bigcup_{n\in\omega}\mbox{{\rm T}}_{N}(\mathcal{K},I,n)$ is equivalent to the set theory 
${\sf ZF}+(V=L)+(\mathcal{K}\mbox{ {\rm is }}  \Pi^{1}_{N}\mbox{{\rm -indescribable}})$.
\end{lemma}

\section{Operator controlled derivations for $\Pi^{1}_{N}$-indescribable cardinals}\label{sect:controlledOme}

In this section,
operator controlled derivations are first introduced,
and
inferences $(\mathcal{K}\in M_{N})$
for $\Pi^{1}_{N}$-indescribability of $\calk$ are then eliminated from operator controlled derivations of $\Pi^{1}_{k}$-sentences $\varphi^{V_{\mathcal{K}}}$ over $\mathcal{K}$.

In what follows $n$ denotes a fixed positive integer.
Let us write
${\sf ZFL}:={\sf ZF}+(V=L)$.

\subsection{Intuitionistic fixed point theories $\mbox{FiX}^{i}(\mbox{{\sf ZFLK}}_{k,n})$}\label{subsec:intfixZFL}

For the fixed positive integer $n$ and $0\leq k< N$,
let
\beqnarr\label{eq:abgam}
&&
b_{n}:=\Psi_{\calk^{+},n}(\ome_{n-1}(I+1)), \:
a_{n}:=\vphi(b_{n})(b_{n}), 
 \nonumber
\\
&&
\gam_{k,n}:=\ome_{n-1}(I+1)+1+a_{n}(N-k),
 \mbox{ and } 
  \nonumber
 \\
 &&
 \fal i<N-k= lh(\vec{\alp}_{k,n})(\vec{\alp}_{k,n}(i)=\gam_{k+i,n})
\eeqnarr
Then ${\sf ZFLK}_{-1,n}:={\sf ZFL}={\sf ZF}+(V=L)$, and for $N>k\geq 0$,
$\mbox{{\sf ZFLK}}_{k,n}, {\sf ZFLK}_{k}$ denotes the set theories 
\beqn\label{eq:ZFLKkn}
{\sf ZFLK}_{k,n}:={\sf ZFL}+(\calk\in Mh_{k,n}(\vec{\alp}_{k,n})),\:
{\sf ZFLK}_{k}:=\bigcup\{{\sf ZFLK}_{k,n}: 0<n<\ome\}
\eeqn
 in the language $\{\in,\mathcal{K},\ome_{1}\}$ with individual constants $\mathcal{K},\ome_{1}$.

Let us also denote the set theory
${\sf ZFL}+(\mathcal{K}\mbox{ is } \Pi^{1}_{N}\mbox{-indescribable})$ in the language $\{\in,\mathcal{K},\ome_{1}\}$
by $\mbox{{\sf ZFLK}}_{N}={\sf ZFLK}_{N,n}$ for any $n$.

To analyze the theory $\mbox{{\sf ZFLK}}_{N}$, we need to handle
a relation
$(\mathcal{H}_{\gamma,n},\Theta,\kap, {\sf ZFLK}_{k,n})\vdash^{a}_{b}\Gamma$ defined in subsection \ref{subsec:operatorcont},
where $n$ is the fixed integer, $k\leq N$,
$\gamma,a,b$ are codes of ordinals with $a<\omega_{n}(I+1)$, $b<I+\omega$,
$\kap$ a regular cardinal and 
$\Theta$ are finite subsets of $L$ and $\Gamma$ a sequent, i.e., a finite set of sentences.
As in \cite{liftupZF, liftupK}  the relation is defined for each $n<\omega$, as a fixed point,
\begin{equation}\label{eq:fixH}
H_{n}(k,\gamma,\Theta,a,b,\kap,\Gamma)\Leftrightarrow
 (\mathcal{H}_{\gamma,n},\Theta,\kap,{\sf ZFLK}_{k,n})\vdash^{a}_{b}\Gamma
\end{equation}

An intuitionistic fixed point theory 
$\mbox{FiX}^{i}(\mbox{{\sf ZFLK}}_{k,n})$ over $\mbox{{\sf ZFLK}}_{k,n}$ is introduced in \cite{intfixset},
and shown to be a conservative extension of $\mbox{{\sf ZFLK}}_{k,n}$.

Fix an $X$-strictly positive formula $\mathcal{Q}(X,x)$ in the language $\{\in,\mathcal{K},\ome_{1},=,X\}$ with an extra unary predicate symbol $X$.
In $\mathcal{Q}(X,x)$ the predicate symbol $X$ occurs only strictly positive.
The language of $\mbox{FiX}^{i}(\mbox{{\sf ZFLK}}_{k,n})$ is $\{\in,\mathcal{K},=,Q\}$ with a fresh unary predicate symbol $Q$.
The axioms in $\mbox{FiX}^{i}(\mbox{{\sf ZFLK}}_{k,n})$ consist of the following:
\benu
\item
All provable sentences in $\mbox{{\sf ZFLK}}_{k,n}$ (in the language $\{\in,\mathcal{K},\ome_{1},=\}$).

\item
Induction schema for any formula $\varphi$ in $\{\in,\mathcal{K},=,Q\}$:
\begin{equation}\label{eq:Qind}
\forall x(\forall y\in x\,\varphi(y)\to\varphi(x))\to\forall x\,\varphi(x)
\end{equation}

\item
Fixed point axiom:
\[
\forall x[Q(x)\leftrightarrow \mathcal{Q}(Q,x)]
.\]
\eenu

The underlying logic in $\mbox{FiX}^{i}(\mbox{{\sf ZFLK}}_{k,n})$ is defined to be the intuitionistic logic.

(\ref{eq:Qind}) yields the following Lemma \ref{lem:vepsfix}.

\begin{lemma}\label{lem:vepsfix}
Let $<^{\varepsilon}$ denote a $\Delta_{1}$-predicate,
which defines a well ordering of type $\veps_{I+1}$.
For each $m<\omega$ and each formula $\varphi$ in $\{\in,\mathcal{K},\ome_{1},=,Q\}$,
\[
\mbox{{\rm FiX}}^{i}(\mbox{{\sf ZFLK}}_{k,n})\vdash\forall x(\forall y<^{\varepsilon}x\,\varphi(y)\to\varphi(x)) \to 
\forall x<^{\varepsilon}\omega_{m}(I+1)\varphi(x)
.\]
\end{lemma}

The following Theorem \ref{th:consvintfix} is shown in \cite{intfixset}.

\begin{theorem}\label{th:consvintfix}
$\mbox{{\rm FiX}}^{i}(\mbox{{\sf ZFLK}}_{k,n})$ is a conservative extension of $\mbox{{\sf ZFLK}}_{k,n}$
for each $k\geq-1$ and each $n\geq 0$.
\end{theorem}

\textit{We will work in} $\mbox{FiX}^{i}(\mbox{{\sf ZFL}})$ \textit{throughout this section with a fixed integer $n$.}

\subsection{Classes of sentences}\label{subsec:classformula}

In this section we consider only the codes of ordinals less than $\omega_{n}(I+1)$ for a fixed positive integer $n$.
$L_{I}:=L=\bigcup_{\alpha\in ORD}L_{\alpha}$ denotes the universe with $\calk\in L$.
$\mbox{{\rm rk}}_{L}(a)=\min\{\alpha\in ORD: a\in L_{\alpha+1}\}$  denotes the $L${\it -rank\/} of $a$.
Also set $\mbox{{\rm rk}}_{L}(L_{I})=I$.

The \textit{language} $\mathcal{L}_{cR}$ is obtained from the language $\{\in, P,P_{I,n},Reg,\mathcal{K},\ome_{1}\}$ 
by adding names (individual constants) $c_{a}$
of each set $a\in L$,
and for each regular cardinal $\kap\leq\calk$ and each subset $B\in\calP(\calk)\cap L_{\calk^{+}}$ 
adding {\it unary predicates\/} $R_{B,\kap}$, their {\it complements\/} $\lnot R_{B,\kap}$
and unary predicate (second-order) variables $X_{i}, \lnot X_{i}\,(i\in\ome)$.
$R_{B,\kap}$ [$\lnot R_{B,\kap}$] denotes the set $B\cap\kap$ [the compliment $\kap\setm B$], resp.
$c_{a}$ is identified with $a$, and
predicate variables are denoted $X,Y,\ldots$
The (individual) variables $x,y,\ldots$ and constants $c_{a}$ are terms.
Terms are denoted $t,s,\ldots$

Then {\it formulae\/} in $\mathcal{L}_{cR}$
are constructed from {\it literals\/} 
$t\in s, t\not\in s, P(t_{1},t_{2},t_{3})$,
\\
$ \lnot P(t_{1},t_{2},t_{3}), P_{I,n}(t), \lnot P_{I,n}(t), Reg(t),\lnot Reg(t), R_{B,\kap}(t),\lnot R_{B,\kap}(t),
X(t),\lnot X(t)
$
by propositional connectives $\lor,\land$, individual quantifiers $\exi x,\fal x$ and 
predicate (second-order) quantifiers $\exi X\subset\kap,\fal X\subset\kap$ for regular cardinals $\kap$ such that
$\ome_{1}<\kap\leq\calk$.
Unbounded quantifiers $\exists x,\forall x$ are denoted by $\exists x\in L_{I},\forall x\in L_{I}$, resp.

For formulae $A$ in $\mathcal{L}_{cR}$, 
${\sf k}(A)$ denotes the set of sets $a\equiv c_{a}$ occurring in $A$, but excluding subsets $B$ in the predicates $R_{B,\kap}$.
${\sf k}(A)\subset L_{I}$ is defined to include 
bounds of `bounded' quantifiers and of `predicates'.
Let us split ${\sf k}(A)$ in two sets ${\sf k}^{E}(A)$ and ${\sf k}^{R}(A)$.
${\sf k}^{R}(A)$
is the set of $\kap$ for which $R_{B,\kap}$ occurs in $A$ for some $B$.
By definition we set $0\in{\sf k}^{R}(A)\cap{\sf k}^{E}(A)$.
In the following definition, $Var$ denotes the set of variables.

\begin{definition}
\benu

\item
${\sf k}^{E}(\lnot A)={\sf k}^{E}(A)$ {\rm and similarly for} ${\sf k}^{R}$.
\item
${\sf k}^{R}(M)=\{0\}$ {\rm for any literal} $M$ {\rm other than} $R_{B,\kap}(t),\lnot R_{B,\kap}(t)$.

\item
${\sf k}^{R}(R_{B,\kap}(t))=\{\kap,0\}$.

\item
${\sf k}^{E}(Q(t_{1},\ldots,t_{m}))=\{t_{1},\ldots,t_{m},0\}\cap L_{I}$
{\rm for literals} $Q(t_{1},\ldots,t_{m})$ {\rm with predicates} $Q$ {\rm in the set} 
$\{P,P_{I,n},Reg,\in\}\cup\{X_{i}\}_{i\in\ome}\cup\{R_{B,\kap}: B\in\calP(\calk)\cap L_{\calk^{+}}, \kap\leq\calk\}$.




\item
${\sf k}^{E}(A_{0}\lor A_{1})={\sf k}^{E}(A_{0})\cup{\sf k}^{E}(A_{1})$ {\rm and similarly for} ${\sf k}^{R}$.

\item
{\rm For} $t\in L_{I}\cup\{L_{I}\}\cup Var$,
${\sf k}^{E}(\exists x\in t\, A(x))=(\{t\}\cap L_{I})\cup{\sf k}^{E}(A(x))$, {\rm while}
${\sf k}^{R}(\exists x\in t\, A(x))={\sf k}^{R}(A(x))$.

\item
${\sf k}^{E}(\exi X\subset \kap\, A(X))=\{\kap\}\cup{\sf k}^{E}(A(X))$, {\rm while}
${\sf k}^{R}(\exi X\subset \kap\, A(X))={\sf k}^{R}(A(X))$.

\item
${\sf k}(A)={\sf k}^{E}(A)\cup{\sf k}^{R}(A)$.
\eenu

\end{definition}

\begin{definition}\label{df:fmlclasses}
\benu

\item
$\Del_{0}$-formulae
{\rm are constructed from literals}
$t\in s, t\not\in s$,
{\rm by propositional connectives} $\lor,\land$, {\rm and} bounded {\rm individual quantifiers} $\exi x\in a,\fal x\in a\,(a\in L_{I})$.
{\rm Note that the predicates} $P,P_{I,n},Reg, R_{B,\kap}, X$
{\rm do not occur in} $\Del_{0}${\rm -formulae.}

\item
{\rm Putting} $\Sigma_{0}:=\Pi_{0}:=\Delta_{0}${\rm , the classes}
$\Sigma_{m}$ {\rm and} $\Pi_{m}$ {\rm of formulae in the language} $\{\in\}$ {\rm with terms}
{\rm are defined as usual using quantifiers}
$\exists x\in L_{I},\forall x\in L_{I}${\rm , where by definition}
$\Sigma_{m}\cup\Pi_{m}\subset\Sigma_{m+1}\cap\Pi_{m+1}$.

{\rm Each formula in} $\Sigma_{m}\cup\Pi_{m}$ {\rm is in prenex normal form with alternating unbounded quantifiers and}
$\Delta_{0}${\rm -matrix.}

{\rm Note that the predicates} $P,P_{I,n},Reg, R_{B,\kap},X$ do not {\rm occur in} $\Sigma_{m}${\rm -formulae.}

\item
{\rm Let} $\lam$ {\rm be a regular cardinal such that} $\ome_{1}<\lam\leq\calk$.

$A\in\Delta_{0}(\lambda)=\Sig_{0}(\lam)=\Pi_{0}(\lam)$ {\rm iff}
{\rm the sentence} $A$ {\rm contains}
{\rm no unbounded quantifiers} $\exi x,\fal x$,
$\max\{\rk_{L}(t): t\in{\sf k}^{R}(A)\}\leq\lam$
 {\rm and}
 $\max\{\rk_{L}(t): t\in {\sf k}^{E}(A)\}<\lambda$,
 {\rm where} $\rk_{L}(x):=0$ {\rm for any variable} $x$.
 
{\rm Note that the predicates}
$P,P_{I},Reg, R_{B,\kap}\, (\kap\leq\lam), X$ {\rm and predicate quantifiers} $\exi X\subset\kap,\fal X\subset\kap\, (\kap<\lam)$ 
{\rm may occur in $\Del_{0}(\lam)$-formulae.}

 

\item 
$A\in \Sigma_{1}(\lambda)$ {\rm iff either} $A\in\Delta_{0}(\lambda)$ {\rm or}
$A\equiv\exists x\in \lambda\, A_{0}$ {\rm with} $A_{0}\in\Delta_{0}(\lambda)$.


\item
{\rm The class of sentences} $\Sigma_{m}(\lambda),\Pi_{m}(\lambda)\, (m<\omega)$ {\rm are defined from}
$\Delta_{0}(\lambda)$ {\rm by applying bounded quantifiers} $\exi x\in\lam,\fal x\in\lam$ {\rm as usual.}

\item
$\Pi^{1}_{0}(\lam):=\bigcup_{m<\ome}\Sig_{m}(\lam)=\bigcup_{m<\ome}\Pi_{m}(\lam)$.

\item 
$A\in \Sigma^{1}_{1}(\lambda)$ {\rm iff either} $A\in\Pi^{1}_{0}(\lambda)$ {\rm or}
$A\equiv\exists X\subset\lambda\, B$ {\rm with} $B\in\Pi^{1}_{0}(\lambda)$.

\item
{\rm The class of sentences} $\Sigma^{1}_{m}(\lambda),\Pi^{1}_{m}(\lambda)\, (m<\omega)$ {\rm are defined from}
$\Pi^{1}_{0}(\lambda)$ {\rm by applying predicate quantifiers} $\exi X\subset\lam,\fal X\subset\lam$ {\rm as usual so that}
$\Sigma^{1}_{m}(\lambda)\cup\Pi^{1}_{m}(\lambda)\subset\Sigma^{1}_{m+1}(\lambda)\cap\Pi^{1}_{m+1}(\lambda)$.

{\rm Set}
$\Pi^{2}_{0}(\lam):=\bigcup_{m<\ome}\Pi^{1}_{m}(\lam)=\bigcup_{m<\ome}\Sig^{1}_{m}(\lam)$.

\item
{\rm (Cf. Definition \ref{df:pi1kfml}.\ref{df:pi1kfml.4}.)}
{\rm For a} $\Pi^{1}_{k}(\lam)${\rm -sentence} $A$ {\rm and a regular cardinal} $\kap$ {\rm with} $\ome_{1}<\kap<\lam$,
$A^{(\kap,\lam)}$ {\rm denotes the result of replacing} 
$R_{B,\lam}(t)$ {\rm by} $R_{B,\kap}(t)$,
{\rm restricting every second-order quantifier} $\exi X\subset \lam, \fal X\subset \lam$ {\rm in} $A$ {\rm to}
$\exi X\subset \kap, \fal X\subset \kap$, {\rm and restricting every bounded quantifier}
$\exists z\in \lam,\forall z\in \lam$ {\rm in} $A$ {\rm to} $\exists z\in \kap, \forall z\in \kap$.

{\rm If} $\max\{\rk_{L}(t): t\in{\sf k}^{R}(A)\}\cap\lam\leq\kap$
 {\rm and}
 $\max\{\rk_{L}(t): t\in {\sf k}^{E}(A)\}<\kap$,
 {\rm then} $A^{(\kap,\lam)}$ {\rm is a} $\Pi^{1}_{k}(\kap)${\rm -sentence.}

{\rm Also for a set} $\Del$ {\rm of} $\Pi^{1}_{k}(\lam)${\rm -sentences},
$\Del^{(\kap,\lam)}:=\{A^{(\kap,\lam)}: A\in\Del\}$.
\eenu

\end{definition}

In what follows we need to consider \textit{sentences} $A$ in the language $\calL_{cR}$.
Sentences are denoted $A,C$ possibly with indices,
while $B,D$ denote sets in $\calP(\calk)\cap L_{\calk^{+}}$.

\begin{definition}
{\rm A set} $\Sigma^{\Sigma_{n+1}}(\lambda)$ {\rm of sentences is defined recursively as follows.}
 \benu

 \item
$\Sigma_{n+1}\subset\Sigma^{\Sigma_{n+1}}(\lambda)$.

 \item
 {\rm Each literal including} $Reg(a), P(a,b,c),P_{I,n}(a), R_{B,\kap}(a)$ {\rm for}
 $\kap<\lam$ 
 {\rm and their negations 
 is in} $\Sigma^{\Sigma_{n+1}}(\lambda)$.
 

 \item
 $\Sigma^{\Sigma_{n+1}}(\lambda)$ {\rm is closed under propositional connectives} $\lor,\land$.

 \item
 {\rm Suppose}
 $\forall x\in b\, A(x)\not\in\Delta_{0}${\rm . Then}
 $\forall x\in b\, A(x)\in \Sigma^{\Sigma_{n+1}}(\lambda)$ {\rm iff} $A(\emptyset)\in \Sigma^{\Sigma_{n+1}}(\lambda)$ {\rm and}
 $\mbox{{\rm rk}}_{L}(b)<\lambda$.

 \item
 {\rm Suppose}
 $\exists x\in b\, A(x)\not\in\Delta_{0}${\rm . Then}
 $\exists x\in b\, A(x)\in \Sigma^{\Sigma_{n+1}}(\lambda)$ {\rm iff} $A(\emptyset)\in \Sigma^{\Sigma_{n+1}}(\lambda)$ {\rm and}
 $\mbox{{\rm rk}}_{L}(b)\leq\lambda$.
 
  \item
 $\exists X\subset\kap\, A(X)\in \Sigma^{\Sigma_{n+1}}(\lambda)$ {\rm iff} 
  $\fal X\subset\kap\, A(X)\in \Sigma^{\Sigma_{n+1}}(\lambda)$ {\rm iff}
  $A(\emptyset)\in \Sigma^{\Sigma_{n+1}}(\lambda)$ {\rm and}
 $\kap<\lambda$.
 
 \eenu
 \end{definition}
 
Note that if
$\kap<\lam$, then any $\Pi^{2}_{0}(\kap)$-sentence is a $\Sig^{\Sig_{n+1}}(\lam)$-formula.


\begin{definition}\label{df:domFfml}
{\rm Let either} $\lam\in Reg$ {\rm and} $x=\Psi_{\lam,n}\bet$, {\rm or} $x=\Psi_{I,n}\bet$ {\rm for some} $\bet$.
{\rm The domain} $dom(F^{\Sigma_{1}}_{x\cup\{\lam\}})=\mbox{{\rm Hull}}_{\Sigma_{1}}^{I}(x\cup\{\lam\})$
{\rm of Mostowski collapse is extended to formulae.}
\[
dom(F^{\Sigma_{1}}_{x\cup\{\lam\}})=\{A\in\Sigma_{1}\cup\Pi_{1}: {\sf k}(A)\subset\mbox{{\rm Hull}}_{\Sigma_{1}}^{I}(x\cup\{\lam\})\}
.\]
{\rm For} $A\in dom(F^{\Sigma_{1}}_{x\cup\{\lam\}})$,
$F^{\Sigma_{1}}_{x\cup\{\lam\}}" A$ {\rm denotes the result of replacing each constant} $\gamma$ {\rm by} 
$F^{\Sigma_{1}}_{x\cup\{\lam\}}(\gamma)${\rm , 
 each unbounded existential quantifier} $\exists z\in L_{I}$ {\rm by} $\exists z\in L_{F^{\Sigma_{1}}_{x\cup\{\lam\}}(I)}${\rm ,
and each unbounded universal quantifier} $\forall z\in L_{I}$ {\rm by} $\forall z\in L_{F^{\Sigma_{1}}_{x\cup\{\lam\}}(I)}$.

{\rm For sequent, i.e., finite set of sentences} $\Gamma\subset dom(F^{\Sigma_{1}}_{x\cup\{\lam\}})${\rm , put}
 $F^{\Sigma_{1}}_{x\cup\{\lam\}}"\Gamma=\{F^{\Sigma_{1}}_{x\cup\{\lam\}}" A: A\in\Gamma\}$.
 
 {\rm Likewise the domain} $dom(F^{\Sigma_{n}}_{x})=\mbox{{\rm Hull}}_{\Sigma_{n}}^{I}(x)$
{\rm is extended to}
\[
dom(F^{\Sigma_{n}}_{x})=\{A\in\Sigma_{n}\cup\Pi_{n}: {\sf k}(A)\subset\mbox{{\rm Hull}}_{\Sigma_{n}}^{I}(x)\}
\]
{\rm and for formula} $A\in dom(F^{\Sigma_{n}}_{x})$,
$F^{\Sigma_{n}}_{x}" A$, {\rm and sequent} $\Gamma\subset dom(F^{\Sigma_{n}}_{x})$,
 $F^{\Sigma_{n}}_{x}"\Gamma$
 {\rm are defined similarly.}
 
\end{definition}

\bprp\label{prp:domFfml}
Suppose that either $\lam\in Reg$ and $x=\Psi_{\lam,n}\bet$, or $x=\Psi_{I,n}\bet$ for some $\bet$.
For $F=F^{\Sig_{1}}_{x\cup\{\lam\}}, F^{\Sig_{n}}_{y}$ and $A\in dom(F)$,
$L_{I}\models A\lrarw F" A$.

\eprp
\bprf
Note that the predicates $P,P_{I,n},Reg, R_{B,\kap},X$ do not occur in $\Sigma_{m}$-formulae.
\eprf
\\

The assignment of disjunctions and conjunctions to sentences is defined as in
\cite{liftupZF} slightly modified, and by adding the clauses for second-order formulae.

Let us define truth values of literals.
\benu
\item
For any literal $M$, $\lnot M$ is \textit{true} iff $M$ is not true.

\item $a\in b$ is \textit{true} iff $L\models a\in b$.

\item
$R_{B,\kap}(a)$  is \textit{true} iff $a\in B\cap\kap$.

\item
 $Reg(a)$ is \textit{true} iff $a\in Reg$, i.e., either $a=\ome_{1}$ or  $\calk<a<I$ is regular.
 
 \item
 
$P(a,b,c)$ is \textit{true} iff 
$a\in Reg$ and $\exists \alpha<\ome_{n}(I+1)[b=\Psi_{a,n}\alpha\,\&\,
c=F^{\Sigma_{1}}_{b\cup\{a\}}(I)]$.

\item
$P_{I,n}(a)$ is \textit{true} iff 
$ \exists \alpha<\ome_{n}(I+1)[a=\Psi_{I,n}\alpha]$.
\eenu

\begin{definition}\label{df:assigndc}
\benu

\item\label{df:assigndc0}
{\rm If} $M$ {\rm is a literal}, {\rm then for} $J:=0$
\[
M:\simeq
\left\{
\begin{array}{ll}
\bigvee(A_{\iota})_{\iota\in J} & \mbox{{\rm if }} M \mbox{ {\rm is false (in }} L_{I}\mbox{{\rm )}}
\\
\bigwedge(A_{\iota})_{\iota\in J} &  \mbox{{\rm if }} M \mbox{ {\rm is true}}
\end{array}
\right.
\]

\item
$(A_{0}\lor A_{1}):\simeq\bigvee(A_{\iota})_{\iota\in J}$
{\rm and}
$(A_{0}\land A_{1}):\simeq\bigwedge(A_{\iota})_{\iota\in J}$
{\rm for} $J:=2$.

\item
{\rm Let} $(\exists z\in b \, \theta[z])\in\Sigma_{n}$ {\rm for} $b\in L_{I}\cup\{L_{I}\}$.
{\rm Then for the set}
\begin{equation}\label{eq:dfmu}
d:=\mu z\in b \, \theta[z] :=
\min_{<_{L}}\{d : (d\in b \land \theta[d]) \lor (\lnot\exists z\in b\, \theta[z]\land d=0)\}
\end{equation}
{\rm with a canonical well ordering} $<_{L}$ {\rm on} $L$,
{\rm and}
$J=\{d\}$
\begin{eqnarray}
\exists z\in b\, \theta[z] & :\simeq & \bigvee(d\in b\land \theta[d])_{d\in J}
\label{eq:sigpimu}
\\
\forall z\in b\, \lnot\theta[z] & :\simeq & \bigwedge(d\in b \to \lnot\theta[d])_{d\in J}
\nonumber
\end{eqnarray}
{\rm where} $d\in b$ {\rm denotes a true literal, e.g.,} $d\not\in d$ {\rm when} $b=L_{I}$.

{\rm This case is applied only when} $\exists z\in b\, \theta[z]$ {\rm is a formula in} $\{\in\}\cup L_{I}${\rm , and}
$(\exists z\in b \, \theta[z])\in\Sigma_{n}$.

\item
{\rm For a sentence} $\exi X\subset\kap\,\tht[X]$ {\rm with a regular}
$\ome_{1}<\kap\leq\calk$ {\rm set}
$
J:=\calP(\calk)\cap L_{\calk^{+}} 
$
{\rm and let}
\[
\exists X\subset\kap\,\tht[X]:\simeq\bigvee(\tht[R_{B,\kap}])_{B\in J}
\mbox{ {\rm and }}
\forall X\subset\kap\,\tht[X]:\simeq\bigwedge(\tht[R_{B,\kap}])_{B\in J}
\]
{\rm where} $\tht[R_{B,\kap}]$ {\rm is obtained from} $\tht[X]$ {\rm by replacing} $X(t)$ $[\lnot X(t)]$
{\rm by} $R_{B,\kap}(t)$ $[\lnot R_{B,\kap}(t)]${\rm , resp.}

\item
{\rm Otherwise set for} $a\in L_{I}\cup\{L_{I}\}$ {\rm and} $J:=\{b: b\in a\}$
\[
\exists x\in a\, A(x):\simeq\bigvee(A(b))_{b\in J}
\mbox{ {\rm and }}
\forall x\in a\, A(x):\simeq\bigwedge(A(b))_{b\in J}
.\]
{\rm This case is applied if one of the predicates} $P,P_{I,n},Reg, R_{B,\kap}$ {\rm occurs in} $\exists x\in a\, A(x)$,
{\rm or} $(\exists x\in a\, A(x))\not\in\Sigma_{n}$.

\eenu

\end{definition}

The definition of the rank $\mbox{{\rm rk}}(A)$ of sentences $A$ 
is modified from \cite{liftupZF, liftupK} so as to the following 
Propositions \ref{lem:rank}.\ref{lem:rank1} and \ref{lem:rank}.\ref{lem:rank3} holds.
The rank $\mbox{{\rm rk}}(A)$ of sentences $A$ is defined by recursion on the  number of symbols occurring in 
$A$.

\begin{definition}\label{df:rank}
\benu

\item\label{df:rank1}
$\mbox{{\rm rk}}(\lnot A):=\mbox{{\rm rk}}(A)$.

\item\label{df:rank2}
$\mbox{{\rm rk}}(M):=\max\{\rk_{L}(t): t\in{\sf k}(M)\}$ {\rm for any (closed) literal} $M$.



\item\label{df:rank6}
$\mbox{{\rm rk}}(A_{0}\lor A_{1}):=\max\{\mbox{{\rm rk}}(A_{0}),\mbox{{\rm rk}}(A_{1})\}+1$.

\item\label{df:rank7}
 \benu
 \item
 $\rk(\exi x<I[b<x\land P_{I,n}(x)]):=I$.
 \item
 $\rk(\exi x<\lam\exi y<\lam[b<x\land P(\lam,x,y)]):=\max\{\lam+1, \rk_{L}(b)\}$ {\rm for} $\lam\in Reg$.
 \item
 $\rk(\exists x\in a\, A(x)):=\max\{\omega \rk_{L}(a), \mbox{{\rm rk}}(A(\emptyset))+1\}$
 {\rm in other cases.}
 \eenu

\item\label{df:rank8}
$
\mbox{{\rm rk}}(\exists X\subset\kap\, A(X)):=
\max\{\kap, \mbox{{\rm rk}}(A(R_{\emptyset,\kap}))+1\}
$.

\eenu

\end{definition}

\begin{proposition}\label{lem:rank}
Let $A$ be a sentence with 
$A\simeq\bigvee(A_{\iota})_{\iota\in J}$ or $A\simeq\bigwedge(A_{\iota})_{\iota\in J}$.
\benu

\item\label{lem:rank15}
$A\in \Sigma^{\Sigma_{n+1}}(\lambda)\Rightarrow\forall\iota\in J(A_{\iota}\in \Sigma^{\Sigma_{n+1}}(\lambda))$.

\item\label{prp:rksig3}
For an ordinal $\lambda\leq I$ with $\omega\lambda=\lambda$,
$\mbox{{\rm rk}}(A)<\lambda \Rightarrow A\in\Sigma^{\Sigma_{n+1}}(\lambda)$.

\item\label{lem:rank0}
$\mbox{{\rm rk}}(A)<I+\omega$.

\item\label{lem:rank1}
$\max\{\rk_{L}(t): t\in{\sf k}(A)\}\leq \mbox{{\rm rk}}(A)\in\{\omega  \mbox{{\rm rk}}_{L}(a)+i  : a\in{\sf k}(A)\cup\{I\}, i\in\omega\}
\subset\mbox{{\rm Hull}}_{\Sigma_{1}}^{I}({\sf k}(A))$.

\item\label{lem:rank2}
$\forall\iota\in J(\mbox{{\rm rk}}(A_{\iota})<\mbox{{\rm rk}}(A))$ if $A$ is not a formula of the form
$\exi x<\lam\exi y<\lam[b<x\land P(\lam,x,y)]$ for some $\lam\in Reg$.

\item\label{lem:rank3}
For an ordinal $\lambda\leq I$ with $\omega\lambda=\lambda$,
$\rk(A)<\lam \Rarw A\in\Del_{0}(\lam)$

\eenu

\end{proposition}

\subsection{Operator controlled derivations}\label{subsec:operatorcont}

\begin{definition}\label{df:operator}
{\rm By an} {\it operator\/} {\rm we mean a map} $\mathcal{H}$, 
$\mathcal{H}:\mathcal{P}(L_{I})\to\mathcal{P}(L_{I}\cup\omega_{n+1}(I+1))${\rm , such that}
\benu
\item
$\forall X\subset L_{I}[X\subset\mathcal{H}(X)]$.

\item
$\forall X,Y\subset L_{I}[Y\subset\mathcal{H}(X) \Rightarrow \mathcal{H}(Y)\subset\mathcal{H}(X)]$.
\eenu

{\rm For an operator} $\mathcal{H}$ {\rm and} $\Theta,\Lambda\subset L_{I}$,
$\mathcal{H}[\Theta](X):=\mathcal{H}(X\cup\Theta)${\rm , and}
$\mathcal{H}[\Theta][\Lambda]:=(\mathcal{H}[\Theta])[\Lambda]${\rm , i.e.,}
$\mathcal{H}[\Theta][\Lambda](X)=\mathcal{H}(X\cup\Theta\cup\Lambda)$.
\end{definition}
Obviously $\mathcal{H}_{\alpha,n}$ is an operator for any $\alp,n$, and
if $\calh$ is an operator, then so is $\mathcal{H}[\Theta]$.
\\

\textit{Sequents} are finite sets of sentences, and inference rules are formulated in one-sided sequent calculus.
Let $\mathcal{H}=\calh_{\gam,n}$ be an operator, $\Tht$ a finite set of subsets of $\calk$,
$\sig\leq I$ a cardinal,
$\Gamma$ a sequent, 
$-2\leq k\leq N$, $a<\omega_{n}(I+1)$ and 
$b<I+\omega$.

We define a relation $(\mathcal{H},\Tht,\sig,{\sf ZFLK}_{k,n})\vdash^{a}_{b}\Gamma$, which is read `there exists an infinitary derivation
of $\Gamma$ which is $(\Tht,\sig,{\sf ZFLK}_{k,n})${\it -controlled\/} by $\mathcal{H}$, and 
 whose height is at most $a$ and its cut rank is less than $b$'.

Recall that $\kap\in Reg$ iff either $\kap>\calk$ is regular or $\kap=\ome_{1}$.
$\kappa,\lambda,\sigma,\pi$ ranges over $Reg^{+}:=Reg\cup\{I\}$.

\begin{definition}\label{df:controlderreg}
$(\mathcal{H},\Tht,\sig,{\sf ZFLK}_{k,n})\vdash^{a}_{b}\Gamma$ {\rm holds if}

\begin{equation}\label{eq:controlder}
{\sf k}(\Gamma)\cup\{a\}\subset\mathcal{H}[\Tht]=\mathcal{H}(\Tht)
\end{equation}

{\rm and one of the following
cases holds:}

\bdes

\item[$(\bigvee)$]
$A\simeq\bigvee\{A_{\iota}: \iota\in J\}$, 
$A\in\Gamma$ {\rm and there exist} $\iota\in J$
{\rm and}
 $a(\iota)<a$ {\rm such that}
\beqn\label{eq:bigveebnd}
\rk_{L}(\iota)<\sig \Rarw \rk_{L}(\iota)< a
\eeqn
{\rm and}
$(\mathcal{H},\Tht,\sig,{\sf ZFLK}_{k,n})\vdash^{a(\iota)}_{b}\Gamma,A_{\iota}$.

\item[$(\bigwedge)$]
$A\simeq\bigwedge\{A_{\iota}: \iota\in J\}$, 
$A\in\Gamma$ {\rm and for every}
$\iota\in J$ {\rm there exists an} $a(\iota)<a$ 
{\rm such that}
$(\mathcal{H},\Tht\cup\{\iota\},\sig,{\sf ZFLK}_{k,n})\vdash^{a(\iota)}_{b}\Gamma,A_{\iota}$.

\item[$(cut)$]
{\rm There exist} $a_{0}<a$ {\rm and} 
$C$
{\rm such that} $\mbox{{\rm rk}}(C)<b$
{\rm and}
$(\mathcal{H},\Tht,\sig,{\sf ZFLK}_{k,n})\vdash^{a_{0}}_{b}\Gamma,\lnot C$
{\rm and}
$(\mathcal{H},\Tht,\sig,{\sf ZFLK}_{k,n})\vdash^{a_{0}}_{b}C,\Gamma$.

\item[$(\mbox{P}_{\lambda})$]
$\lambda\in Reg$ {\rm and there exists} $\alpha<\lambda$ {\rm such that} 
$(\exists x<\lambda\exists y<\lambda[\alpha<x \land P(\lambda,x,y)])\in\Gamma$.

\item[$(\mbox{F}^{\Sigma_{1}}_{x\cup\{\lambda\}})$] 
$\lambda\in Reg$,
$x=\Psi_{\lambda,n}\beta\in\mathcal{H}$ {\rm for a} $\beta$
{\rm and there exist} $a_{0}<a$, 
$\Gamma_{0}\subset\Sigma_{1}$ {\rm and} $\Lambda$ {\rm such that}
${\sf k}(\Gamma_{0})\subset\mbox{{\rm Hull}}_{\Sigma_{1}}^{I}((\mathcal{H}\cap x)\cup\{\lambda\})$,
$\Gamma=\Lambda\cup (F^{\Sigma_{1}}_{x\cup\{\lambda\}}"\Gamma_{0})$
{\rm and}
$
(\mathcal{H},\Tht,\sig,{\sf ZFLK}_{k,n})\vdash^{a_{0}}_{b}\Lambda,\Gamma_{0}$,
{\rm where} $F^{\Sigma_{1}}_{x\cup\{\lambda\}}$ {\rm denotes the Mostowski collapse}
$F^{\Sigma_{1}}_{x\cup\{\lambda\}}: \mbox{{\rm Hull}}^{I}_{\Sigma_{1}}(x\cup\{\lambda\})\leftrightarrow L_{F^{\Sigma_{1}}_{x\cup\{\lambda\}}(I)}$.

\item[$(\mbox{P}_{I,n})$]
{\rm There exists} $\alpha<I$ {\rm such that} 
$(\exists x<I[\alpha<x \land P_{I,n}(x)])\in\Gamma$.

\item[$(\mbox{F}^{\Sigma_{n}}_{x})$] 
$
x=\Psi_{I,n}\beta\in\mathcal{H}$
{\rm for a} $\beta$ {\rm and there exist} $a_{0}<a$, $\Gamma_{0}\subset\Sigma_{n}$
{\rm and} $\Lambda$ {\rm such that}
${\sf k}(\Gamma_{0})\subset\mbox{{\rm Hull}}_{\Sigma_{n}}^{I}(\mathcal{H}\cap x)$,
$\Gamma=\Lambda\cup (F^{\Sigma_{n}}_{x}"\Gamma_{0})$ {\rm and}
$(\mathcal{H},\Tht,\sig,{\sf ZFLK}_{k,n})\vdash^{a_{0}}_{b}\Lambda,\Gamma_{0}$,
{\rm where} $F^{\Sigma_{n}}_{x}$ {\rm denotes the Mostowski collapse}
$F^{\Sigma_{n}}_{x}: \mbox{{\rm Hull}}^{I}_{\Sigma_{n}}(x)\leftrightarrow L_{F^{\Sigma_{n}}_{x}(I)}$.

\item[$(\exi^{2}(\kap))$]$(\ome_{1}<\kap\leq\calk)$

$\exi X\subset\kap\,\tht[X]\simeq\bigvee(\tht[R_{B,\kap}])_{B\in J}$, $(\exi X\subset\kap\,\tht[X])\in\Gam\cap\Pi^{2}_{0}(\kap)$ {\rm and there exist}
$B\in J=\calP(\calk)\cap L_{\calk^{+}}$
{\rm and} $a_{0}<a$ {\rm such that} 
$
(\calh,\Tht,\sig,{\sf ZFLK}_{k,n})\vdash^{a_{0}}_{b}\Gam,\tht[R_{B,\kap}]
$.

\item[$(\fal^{2}(\kap))$]$(\ome_{1}<\kap\leq\calk)$

$\fal X\subset\kap\,\tht[X]\simeq\bigwedge(\tht[R_{B,\kap}])_{B\in J}$, $(\fal X\subset\kap\,\tht[X])\in\Gam\cap\Pi^{2}_{0}(\kap)$ {\rm and}
{\rm there exists an} $a_{0}<a$ {\rm such that}
{\rm for any} $B\in J=\calP(\calk)\cap L_{\calk^{+}}$
$(\calh,\Tht,\sig,{\sf ZFLK}_{k,n})\vdash^{a_{0}}_{b}\Gam,\tht[R_{B,\kap}]$.

{\rm Note that the ordinal} $a_{0}$ {\rm is independent from} $B\in J$.






\item[$(\pi\in Mh_{i,n}(\vec{\alp})\mbox{{\rm [}}\Tht \mbox{{\rm ]}},\vec{\nu})$]
This inference rule is only for $k\geq -1$.

{\rm There exist a regular cardinal} $\pi<\calk$, {\rm a number} $i< N$, {\rm sequences of ordinals} 
$\vec{\alp}=(\alp_{0},\ldots,\alp_{N-i-1}), \vec{\nu}=(\nu_{0},\ldots,\nu_{N-i-1})$ {\rm such that}
$lh(\vec{\alp})=lh(\vec{\nu})=N-i$.
{\rm Moreover there are ordinals} $a_{\ell},a_{r}(\rho), a_{0}$,  
{\rm and a finite set} $\Del$ {\rm of} $\Sig^{1}_{i}(\pi)${\rm -sentences enjoying the following conditions:}

\benu
\item
$\vec{\nu}\in\calh_{\vec{\nu},n}[\Tht]\cap\vec{\alp}(:\Lrarw\fal j<lh(\vec{\nu})(\nu_{j}\in\calh_{\nu_{j},n}[\Tht]\cap\alp_{j}))$,
$\pi\in Mh_{i,n}(\vec{\alp})[\Tht]$, {\rm and}
$b\geq\pi$.

 \item
 \beqn\label{eq:coeff}
 \pi\in\calh[\Tht] 
\eeqn

 \item
 \beqn\label{eq:qbnd}
\fal a\in K(\vec{\alp})(a\leq\gam) 
 \eeqn
 {\rm where} $\calh=\calh_{\gam,n}$.

 \item
 {\rm For each} $\del\in\Del$
 \[
(\mathcal{H},\Tht,\sig,{\sf ZFLK}_{k,n})\vdash^{a_{\ell}}_{b}\Gamma, \lnot\del
\]

\item

{\rm Let} $H_{i}(\vec{\nu},\gam)$ {\rm denote the} \textit{resolvent class} {\rm for}  
$\pi\in Mh_{i,n}(\vec{\alp})\mbox{{\rm [}}\Tht \mbox{{\rm ]}}$
{\rm with respect to} $\vec{\nu}$ {\rm and}  $\gam${\rm :}
\beqn
\label{eq:Hdfk}
H_{i}(\vec{\nu},\gam)  := \{\rho\in \bigcap_{j<lh(\vec{\alp})}Mh_{i+j,n}(((\vec{\nu}\bullet\vec{\alp})[j] )[\Tht]\cap\pi 
:
\calh_{\gam,n}[\Tht](\rho)\cap\pi\subset\rho\} 
\eeqn
{\rm where} $\calh=\calh_{\gam,n}$.
{\rm Then for any} $\rho\in H_{i}(\vec{\nu},\gam)$
\[
(\mathcal{H},\Tht\cup\{\rho\},\sig,{\sf ZFLK}_{k,n})\vdash^{a_{r}(\rho)}_{b}\Gamma, 
\Del^{(\rho,\pi)}
\]

\item
\beqn\label{eq:Mhordgap}
\sup\{a_{\ell},a_{r}(\rho): \rho\in H_{k}(\vec{\nu},\gam)\}\leq a_{0}\in\calh[\Tht]\cap a
\eeqn

\eenu

\item[$(\calk\in Mh_{k,n}(\vec{\alp})\mbox{{\rm [}}\Tht \mbox{{\rm ]}},\vec{\nu})$]
This inference rule is only for $k\geq 0$.

{\rm There are sequences of ordinals} 
$\vec{\alp}=(\alp_{0},\ldots,\alp_{N-k-1}), \vec{\nu}=(\nu_{0},\ldots,\nu_{N-k-1})$ {\rm such that}
$lh(\vec{\alp})=lh(\vec{\nu})=N-k$.
{\rm Moreover there are ordinals} $a_{\ell},a_{r}(\rho), a_{0}$,  
{\rm and a finite set} $\Del$ {\rm of} $\Sig^{1}_{k}(\calk)${\rm -sentences enjoying the following conditions:}

\benu
\item
$\vec{\nu}\in\calh_{\vec{\nu},n}[\Tht]\cap\vec{\alp}$ {\rm and} $b\geq\calk$.

{\rm Note that} $\calk\in Mh_{k,n}(\vec{\alp})[\Tht]$ {\rm is} not required {\rm here.}

 \item
 \beqn
  \renewcommand{\theequation}{\ref{eq:qbnd}}
\fal i< lh(\vec{\alp})(\vec{\alp}(i)\leq\min\{\gam_{k+i,n},\gam\})
 \eeqn
 \addtocounter{equation}{-1}
 {\rm where} $\calh=\calh_{\gam,n}$ {\rm and}
 $\vec{\alp}_{k,n}(i)=\gam_{k+i,n}$,  {\rm cf. (\ref{eq:abgam}) for the ordinal} $\gam_{k,n}$.

 \item
 {\rm For each} $\del\in\Del$
 \[
(\mathcal{H},\Tht,\sig,{\sf ZFLK}_{k,n})\vdash^{a_{\ell}}_{b}\Gamma, \lnot\del
\]

\item

{\rm Let} $H_{k}(\vec{\nu},\gam)$ {\rm denote the} \textit{resolvent class} {\rm for}  
$\calk\in Mh_{k,n}(\vec{\alp})\mbox{{\rm [}}\Tht \mbox{{\rm ]}}$
{\rm with respect to} $\vec{\nu}$ {\rm and}  $\gam${\rm :}
{\small
 \beqn
  \renewcommand{\theequation}{\ref{eq:Hdfk}}
H_{k}(\vec{\nu},\gam)  := \{\rho\in \bigcap_{i<lh(\vec{\alp})}Mh_{k+i,n}(((\vec{\nu}\bullet\vec{\alp})[i] )[\Tht]\cap\calk
:
\calh_{\gam,n}[\Tht](\rho)\cap\calk\subset\rho\} 
 \eeqn
 \addtocounter{equation}{-1}
}
{\rm where} $\calh=\calh_{\gam,n}$.

{\rm Then for any} $\rho\in H_{k}(\vec{\nu},\gam)$
\[
(\mathcal{H},\Tht\cup\{\rho\},\sig,{\sf ZFKL}_{k,n})\vdash^{a_{r}(\rho)}_{b}\Gamma, 
\Del^{(\rho,\calk)}
\]

{\rm In particle when} $k=N$,

\[
H_{N}(\gam):=H_{N}(\vec{\nu},\gam)  = \{\rho<\calk
:
\rho \mbox{ {\rm is regular and }} \calh_{\gam,n}[\Tht](\rho)\cap\calk\subset\rho\} 
\]

\item
\beqn
 \renewcommand{\theequation}{\ref{eq:Mhordgap}}
\sup\{a_{\ell},a_{r}(\rho): \rho\in H_{k}(\vec{\nu},\gam)\}\leq a_{0}\in\calh[\Tht]\cap a
 \eeqn
 \addtocounter{equation}{-1}

\eenu

\edes

\end{definition}
Some comments on  the inference rule $(\pi\in Mh_{i,n}(\vec{\alp})[\Tht],\vec{\nu})$ are helpful.

{\rm The inference rule} $(\pi\in Mh_{i,n}(\vec{\alp})[\Tht],\vec{\nu})$ {\rm says that} $H_{i}(\vec{\nu},\gam)$ 
{\rm is} $\Pi^{1}_{i}${\rm -indescribable in} $\pi$, {\rm and is depicted as follows by suppressing} $\sig,{\sf ZFLK}_{k,n}${\rm :}
{\footnotesize
\[
\infer[(\pi\in Mh_{i,n}(\vec{\alp})\mbox{$[\Tht]$},\vec{\nu})]
{(\mathcal{H},\Tht)\vdash^{a}_{b}\Gam}
{
\{(\mathcal{H},\Theta)\vdash^{a_{\ell}}_{b}\Gamma,\lnot\del\}_{\del\in\Del}
&
 \{(\mathcal{H},\Theta\cup\{\rho\})\vdash^{a_{r}(\rho)}_{b}\Gamma, \Del^{(\rho,\pi)}\}_{
 \rho\in H_{i}(\vec{\nu},\gam)
 \cap \pi}
}
\]
}

{\rm The inference rule}
 $(\calk\in Mh_{N,n}(\vec{\alp})\mbox{{\rm [}}\Tht \mbox{{\rm ]}})$ {\rm is denoted by}
$(\calk\in M_{N})$, {\rm and depicted as follows:}
{\small
\[
\infer[(\calk\in M_{N})]{(\mathcal{H},\Tht)\vdash^{a}_{b}\Gam}
{
\{(\mathcal{H},\Tht)\vdash^{a_{\ell}}_{b}\Gamma, \lnot\del\}_{\del\in\Del}
&
\{(\mathcal{H},\Tht\cup\{\rho\})\vdash^{a_{r}(\rho)}_{b}\Gamma, 
\Del^{(\rho,\calk)}\}_{\rho\in H_{N}(\gam)}
}
\]
}
{\rm where} $\Del$ {\rm is a finite set of} $\Sig^{1}_{N}${\rm -sentences.}

Note that there occurs only the inference rules 
$(\calk\in Mh_{k,n}(\vec{\alp})\mbox{{\rm [}}\Tht \mbox{{\rm ]}},\vec{\nu})$
in the derivation establishing the fact
$(\calh,\Tht,\kap,{\sf ZFLK}_{k,n})\vdash^{a}_{b}\Gam$ for $k\geq 0$.
In particular there occurs no such inference rules for any $k\geq 0$
in the derivation establishing the fact
$(\calh,\Tht,\kap,{\sf ZFLK}_{-1,n})\vdash^{a}_{b}\Gam$.

When an inference rule $(\pi\in Mh_{i,n}(\vec{\alp})\mbox{{\rm [}}\Tht \mbox{{\rm ]}},\vec{\nu})$ is applied for $\pi<\calk$,
$\pi\in Mh_{i,n}(\vec{\alp})[\Tht]$ has to be enjoyed.
As contrasted with this case $\pi<\calk$, 
we don't assume that $\calk\in Mh_{k,n}(\vec{\alp})[\Tht]$ holds in applying inference rules
 $(\calk\in Mh_{k,n}(\vec{\alp})[\Tht],\vec{\nu})$.
 Therefore each inference rule in the derivation establishing the fact
$(\calh,\Tht,\kap,{\sf ZFLK}_{-1,n})\vdash^{a}_{b}\Gam$ is \textit{correct} demonstrably in {\sf ZFL}
in the sense that if (the disjunction of) each upper sequent holds, then so is its lower sequent.
Moreover for $k\geq 0$ and $\Tht\subset(\calk+1)$
each inference rule $(\calk\in Mh_{k,n}(\vec{\alp})[\Tht],\vec{\nu})$
 is correct demonstrably in ${\sf ZFLK}_{k,n}$, cf. Proposition \ref{prp:clshull}.\ref{prp:L4.10.1.1c} and (\ref{eq:ZFLKkn}).
\\

\noindent
An inspection to Definition \ref{df:controlderreg}
shows that there exists a strictly positive formula $H_{n}$
such that the relation $(\mathcal{H}_{\gamma,n},\Theta,\kap,{\sf ZFLK}_{k,n})\vdash^{a}_{b}\Gamma$
is a fixed point of $H_{n}$ as in (\ref{eq:fixH}).

In what follows the relation should be understood as a fixed point of $H_{n}$,
and recall that we are working in 
the intuitionistic fixed point theory $\mbox{FiX}^{i}({\sf ZFL})$
over ${\sf ZFL}$ defined in subsection \ref{subsec:intfixZFL}.


We will state some lemmata for the operator controlled derivations with sketches of their proofs
since
these can be shown as in \cite{Buchholz, liftupZF}.

In what follows by an operator $\mathcal{H}$ we mean an $\mathcal{H}_{\gamma,n}$ for an ordinal $\gamma$.
Also except otherwise stated, $(\calh,\Tht)\vdash^{a}_{b}\Gam$ [$(\calh,\Tht,\kap)\vdash^{a}_{b}\Gam$] 
denotes $(\calh,\Tht,\kap,{\sf ZFLK}_{k,n})\vdash^{a}_{b}\Gam$
for some arbitrarily fixed $\kap$ and $k,n$ [for some fixed arbitrarily $k,n$], resp.

\blem\label{lem:2ndinversion}{\rm (Inversion lemma for predicate quantifiers)}\\
Let $(\fal X\subset\pi\, \tht[X])\in\Pi^{2}_{0}(\pi)$ and assume
$(\calh,\Tht)\vdash^{a}_{b}\Gam,\fal X\subset\pi\, \tht[X]$.
Then for any subset $B\in\calP(\calk)\cap L_{\calk^{+}}$,
$(\calh,\Tht)\vdash^{a}_{b}\Gam,\tht[R_{B,\pi}]$.
\elem



\begin{lemma}\label{lem:tautology}{\rm (Tautology)}
$(\mathcal{H},{\sf k}(\Gamma\cup\{A\}),I)\vdash^{2\footnotesize{\mbox{{\rm rk}}}(A)}_{0}\Gamma,\lnot A, A$.
\end{lemma}




To interpret the axiom of $\Pi^{1}_{N}$-indescribability of $\calk$
\begin{equation}
\renewcommand{\theequation}{\ref{eq:Kord}}
\fal x\subset\calk[\vphi(x)\to\exi \rho<\calk(
\vphi^{(\rho,\calk)}(x\cap\rho))]
\end{equation}
\addtocounter{equation}{-1}
we need to show the equivalence of $\Pi^{1}_{k}(\kap)$-formula and its translation in the first-order language $\{\in\}$.

Let $\vphi(X)$ be a $\Pi^{1}_{N}$-sentence possibly with a predicate constant $X$ in the language $\{\in\}\cup\{X_{i}\}_{i<\ome}$, 
$\vphi_{0}(x)$ its $\Pi^{1}_{N}(\calk)$-translation in $\{\in\}$ defined in Definition \ref{df:pi1kfml}.\ref{df:pi1kfml.3}, and
$\vphi_{1}(X)$ denote the $\Pi^{1}_{N}(\calk)$-sentence obtained from $\vphi(X)$ by
restricting second-order quantifiers $\exi Y,\fal Y$ to $\exi Y\subset\calk, \fal Y\subset\calk$,
and restricting first-order (unbounded) quantifiers $\exi z,\fal z$ to $\exi z\in\calk,\fal z\in\calk$.

Let us temporarily introduce a complexity measure $d(\vphi)<\ome$ of second-order formulae $\vphi$ in $\{\in\}\cup\{X_{i}\}_{i}$.
$d(x\in y)=d(X(y))=0$, $d(\vphi_{0}\lor\vphi_{1})=\max\{d(\vphi_{i}) : i<2\}+1$,
$d(\exi x\,\vphi)=d(\vphi)+1$, $d(\exi X\, \vphi)=d(\vphi)+2$ and $d(\lnot\vphi)=d(\vphi)$.

\bprp\label{prp:predset}
For any set $B\in\calP(\calk)\cap L_{\calk^{+}}$, any $\rho\leq\calk$ and any $\Gam$
\[
(\calh,{\sf k}(\Gam)\cup\{B,\rho\},I)\vdash^{2 d(\vphi)}_{0}\Gam,\lnot\vphi_{0}^{(\rho,\calk)}(B\cap\rho), \vphi_{1}^{(\rho,\calk)}(R_{B,\rho})
\]
and
\[
(\calh,{\sf k}(\Gam)\cup\{B,\rho\},I)\vdash^{2 d(\vphi)}_{0}\Gam,\vphi_{0}^{(\rho,\calk)}(B\cap\rho), \lnot\vphi_{1}^{(\rho,\calk)}(R_{B,\rho})
\]
\eprp
\bprf
This is seen by induction on $d(\vphi)$. 
Let us check one half of the case when $\vphi(X)\equiv(\exi Y\,\tht(Y,X))$.
Let $\Tht={\sf k}(\Gam)\cup\{B,\rho\}$ and $d=d(\tht)$.
By IH we have for any $D\subset\calk$
\[
(\calh,\Tht\cup\{D\},I)\vdash^{2d}_{0}\Gam,\lnot\tht_{0}^{(\rho,\calk)}(D\cap\rho,B\cap\rho), \tht_{1}^{(\rho,\calk)}(R_{D,\rho},R_{B,\rho})
\]
Hence for any $D\subset\rho$
\[
(\calh,\Tht\cup\{D\},I)\vdash^{2d+2}_{0}\Gam,(D\not\subset\rho)\lor \lnot\tht_{0}^{(\rho,\calk)}(D,B\cap\rho), 
\exi Y\subset\rho\, \tht_{1}^{(\rho,\calk)}(Y,R_{B,\rho})
\]
On the other side, if $a\not\subset\rho$, then 
$(\calh,\Tht\cup\{a\},I)\vdash^{1}_{0}a\not\subset\rho$.
Hence
\[
(\calh,\Tht\cup\{a\},I)\vdash^{2}_{0}(a\not\subset\rho)\lor \lnot\tht_{0}^{(\rho,\calk)}(a,B\cap\rho) 
\]
By a $(\bigwedge)$ for $\lnot\vphi_{0}^{(\rho,\calk)}(B\cap\rho)\equiv\fal y(y\subset\rho\to \lnot \tht_{0}^{(\rho,\calk)}(y,B\cap\rho))$
we obtain
\[
(\calh,\Tht,I)\vdash^{2d+3}_{0}\Gam,\fal y(y\subset\rho\to \lnot\tht_{0}^{(\rho,\calk)}(y,B\cap\rho)), 
\exi Y\subset\rho\, \tht_{1}^{(\rho,\calk)}(Y,R_{B,\rho})
\]
\eprf

\begin{definition}
{\rm For a formula} $\exists x\in d\, A$ {\rm and ordinals} $\lambda=\mbox{{\rm rk}}_{L}(d)\in Reg\cup\{I\}, \alpha$,
$(\exists x\in d\, A)^{(\exists\lambda\!\upharpoonright\!\alpha)}$ {\rm denotes the result of restricting the} {\it outermost existential quantifier\/}
$\exists x\in d$ {\rm to} $\exists x\in L_{\alpha}$,
$(\exists x\in d\, A)^{(\exists\lambda\!\upharpoonright\!\alpha)}\equiv
(\exists x\in L_{\alpha}\, A)$.
\end{definition}

In what follows $F_{x,\lambda}$ denotes $F^{\Sigma_{1}}_{x,\lambda}$ when $\lambda\in R$, and $F^{\Sigma_{n}}_{x}$ when $\lambda=I$.

\begin{lemma}\label{lem:boundednessreg}{\rm (Boundedness)}
\\
Let $\lambda\in Reg\cup\{I\}$, $C\equiv (\exists x\in d\, A)$ and $C\not\in\{\exists x<\lambda \exists y<\lambda[\alpha<x \land P(\lambda,x,y)]: \alpha<\lambda\in Reg\}\cup\{\exists x<I[\alpha<x \land P_{I,n}(x)]:\alpha<I\}$.
Assume that $\mbox{{\rm rk}}(C)=\lambda=\mbox{{\rm rk}}_{L}(d)$ and $C$ is not a second-order formula.

\benu

\item\label{lem:boundednessregexi}

\[
(\mathcal{H},\Tht,\lambda)\vdash^{a}_{c}\Lambda, C \,\&\, a\leq b\in\mathcal{H}\cap\lambda
\Rightarrow (\mathcal{H},\Tht,\lambda)\vdash^{a}_{c}\Lambda,C^{(\exists\lambda\!\upharpoonright\! b)}
.\]
\item\label{lem:boundednessregfal}
\[
(\mathcal{H},\Tht,\lam)\vdash^{a}_{c}\Lambda,\lnot C \,\&\, b\in\mathcal{H}\cap\lambda 
\Rightarrow (\mathcal{H},\Tht,\lam)\vdash^{a}_{c}\Lambda,\lnot (C^{(\exists\lambda\!\upharpoonright\! b)})
.\]
\eenu

\end{lemma}





In the following Lemma \ref{lem:predcereg}, note that 
$\rk(\exi x<\lam\exi y<\lam[\alp<x \land P(\lam,x,y)])=\lam+1$ for $\alp<\lam\in Reg$, and
$\rk(\exi x<I[\alp<x \land P_{I,n}(x)])=I$.

\begin{lemma}\label{lem:predcereg}{\rm (Predicative Cut-elimination)}
\benu

\item\label{lem:predcereg2}
Suppose 
$[c,c+\omega^{a}[\cap (\{\lambda+1:\lam\in Reg\}\cup\{I\})=\emptyset$
and $k\neq -2 \Rarw ]c,c+\ome^{a}]\cap\{\kap\leq\calk:\ome_{1}<\kap \mbox{ {\rm is regular}}\}=\emptyset$. Then
$(\mathcal{H},\Tht,\kap,{\sf ZFLK}_{k,n})\vdash^{b}_{c+\omega^{a}}\Gamma
\,\&\, a\in\mathcal{H}[\Tht]
\Rightarrow (\mathcal{H},\Tht,\kap,{\sf ZFLK}_{k,n})\vdash^{\varphi ab}_{c}\Gamma$.

\item\label{lem:predcereg4}
For $\lambda\in Reg$,
$(\mathcal{H}_{\gamma},\Tht,\kap,{\sf ZFLK}_{k,n})\vdash^{b}_{\lambda+2}\Gamma \,\&\, \gamma\in\mathcal{H}_{\gamma}\,\&\,
 \Rightarrow 
(\mathcal{H}_{\gamma+b},\Tht,\kap,{\sf ZFLK}_{k,n})\vdash^{\omega^{b}}_{\lambda+1}\Gamma$.

\item\label{lem:predcereg5}
$(\mathcal{H}_{\gamma},\Tht,\kap,{\sf ZFLK}_{k,n})\vdash^{b}_{I+1}\Gamma \,\&\, \gamma\in\mathcal{H}_{\gamma}\,\&\,
 \Rightarrow 
(\mathcal{H}_{\gamma+b},\Tht,\kap,{\sf ZFLK}_{k,n})\vdash^{\omega^{b}}_{I}\Gamma$.

\eenu

\end{lemma}

The following  Lemma \ref{th:Collapsingthmreg1K} is seen as in  Lemma 4.10 of \cite{liftupK}.

\blem\label{th:Collapsingthmreg1K}{\rm (Collapsing)}\\
Assume 
$\calk<\lam\leq\sig\in\{\ome_{\alp}: \alp\leq I\}\spand\lam\in Reg$, or
$\ome_{1}=\lam\leq\sig\in\{\ome_{\alp}: \alp\leq I\}$ and $k=-2$.

Suppose $\{\gam,\lam,\sig\}\subset\calh_{\gam,n}[\Tht]$ with $\fal\rho\geq\lam[\Tht\subset\calh_{\gam,n}(\Psi_{\rho,n}\gam)]$, 
and 
$
\Gam\subset\Sig^{\Sig_{n+1}}(\lam)
$.
Let
$\mu=\sig+1$ if either $\sig$ is a regular cardinal or $\sig=I$.
Otherwise $\mu=\sig$.

Moreover assume
\[
 (\calh_{\gam,n},\Tht,\sig,{\sf ZFLK}_{k,n})\vdash^{a}_{\mu}\Gam
 \]

Then for $\hat{a}=\gam+\ome^{\sig(1+a)}$ 
we have
\[
  (\calh_{\hat{a}+1,n}, \Tht,\lam,{\sf ZFLK}_{k,n})\vdash^{\Psi_{\lam,n}\hat{a}}_{\Psi_{\lam,n}\hat{a}}\Gam.
\]
\elem

\subsection{Lowering and eliminating higher Mahlo operations}\label{subsec:elimpi11}

In the section we eliminate inferences $(\mathcal{K}\in M_{N})$
for $\Pi^{1}_{N}$-indescribability.

In the following Lemma \ref{lem:CollapsingthmKR}, let for the fixed $n$
\[
(\mathcal{H},\Tht,{\sf ZFLK}_{k,n})\vdash^{a}_{b}\Gamma
 :\Leftrightarrow 
(\mathcal{H},\Tht,\calk^{+},{\sf ZFLK}_{k,n})\vdash^{a}_{b}\Gamma.
\]

Recall that we have defined ordinals in (\ref{eq:abgam}) as follows:
$b_{n}:=\Psi_{\calk^{+},n}(\ome_{n-1}(I+1))$, $a_{n}:=\vphi(b_{n})(b_{n})$, $\gam_{k,n}:=\ome_{n-1}(I+1)+1+a_{n}(N-k)$,
and $\fal i<N-k= lh(\vec{\alp}_{k,n})(\vec{\alp}_{k,n}(i)=\gam_{k+i,n})$.

\begin{lemma}\label{lem:CollapsingthmKR}{\rm ($\mbox{FiX}^{i}({\sf ZFL})$)}\\
Let $k>0$.
Suppose for an operator $\mathcal{H}_{\gam,n}$ and an ordinal $a\leq a_{n}$ 
\beqn\label{eq:CollapsingthmKR00}
(\mathcal{H}_{\gam,n},\Tht,{\sf ZFLK}_{k,n})\vdash^{a}_{\calk}\Lam,\Gamma
\eeqn
where $\gam\in\calh_{\gam,n}[\Tht]$ and $\gam\leq\gam_{k,n}$,
$\Lam\subset\Sig^{\Sig_{n}}(\lam)$ for some regular cardinal $\lam\in\calh_{\gam,n}[\Tht]\cap\calk$, 
$\Gamma$ consists of $\Pi^{1}_{k}(\calk)$-sentences.

Let
\[
\hat{a} := \gam+a\leq\gam_{k-1,n}
\]
Then 
for any $\kap$ such that either $\kap=\calk$, or 
$\kap\in Mh_{k-1,n}((\hat{a})*\vec{\alp}_{k,n})[\Tht\cup\{\kap\}]\cap\calk$ and
$\calh_{\hat{a},n}[\Tht](\kap)\cap\calk\subset\kap$,

\beqn\label{eq:CollapsingthmKR01}
(\mathcal{H}_{\hat{a}},\Tht\cup\{\kap\},{\sf ZFLK}_{k-1,n})\vdash^{\kap+\ome a}_{\kap}
\Lam,\Gamma^{(\kap,\calk)}
\eeqn
holds.

\end{lemma}
{\bf Proof} by induction on $a$. 
Let $\mathcal{H}=\mathcal{H}_{\gam,n}$.

Let either $\kap=\calk$, or $\kap\in Mh_{k-1,n}((\hat{a})*\vec{\alp}_{k,n})[\Tht\cup\{\kap\}]$, $\kap<\calk$
and
$\calh_{\hat{a},n}[\Tht](\kap)\cap\calk\subset\kap$.
From $\calh_{\hat{a},n}[\Tht](\kap)\cap\calk\subset\kap$ and $\gam\leq\hat{a}$ we see that
\beqn\label{eq:CollapsingthmKR100}
\lam\in\calh[\Tht](\kap)\cap\calk\subset \kap
\eeqn
Also by (\ref{eq:controlder}) we have ${\sf k}(\Gam)\subset\calh[\Tht]\cap(\calk+1)$ for the set $\Gam$ of $\Pi^{1}_{k}(\calk)$-sentences.
Hence $\Gam^{(\kap,\calk)}$ is a set of $\Pi^{1}_{k}(\kap)$-sentences.

By $\hat{a}\geq\gam$, we obtain $\fal b\in K((\hat{a})*\vec{\alp})(b\leq\hat{a})$ if 
$\fal b\in K(\vec{\alp})(b\leq\gam)$, cf.  (\ref{eq:qbnd}).
\\

\noindent
{\bf Case 1}. 
First consider the case when the last inference is a $(\calk\in Mh_{k,n}(\vec{\alp})[\Tht],\vec{\nu})$
where $\fal i<N-k(\vec{\alp}(i)\leq\vec{\alp}_{k,n})$.
We have $a_{\ell}\in\mathcal{H}[\Tht]\cap a$, and
$a_{r}(\rho)\in\mathcal{H}[\Tht\cup\{\rho\}]\cap a$.
$\Del$ is a finite set of $\Sig^{1}_{k}(\calk)$-sentences.

{\small
\[
\infer{(\mathcal{H},\Tht,{\sf ZFLK}_{k,n})\vdash^{a}_{\calk}\Lam,\Gam}
{
\{(\mathcal{H},\Theta,{\sf ZFLK}_{k,n})\vdash^{a_{\ell}}_{\calk}\Lam,\Gamma,\lnot\del\}_{\del\in\Del}
&
 \{(\mathcal{H},\Theta\cup\{\rho\},{\sf ZFLK}_{k,n})\vdash^{a_{r}(\rho)}_{\calk}
 \Lam,\Gamma, \Del^{(\rho,\calk)}\}_{ \rho\in H_{k}(\vec{\nu},\gam)}
}
\]
}
\noindent
where $H_{k}(\vec{\nu},\gam)$ is the resolvent class for $\calk\in Mh_{k,n}(\vec{\alp})[\Tht]$ with respect to a 
$\vec{\nu}$ and $\gam$:
\beqn
\renewcommand{\theequation}{\ref{eq:Hdfk}}
H_{k}(\vec{\nu},\gam)  := 
\{\rho\in \bigcap_{i<lh(\vec{\alp})}Mh_{k+i,n}((\vec{\nu}\bullet\vec{\alp})[i] )[\Tht]\cap\calk
:
\calh_{\gam,n}[\Tht](\rho)\cap\calk\subset\rho 
\}
\eeqn
\addtocounter{equation}{-1}
and $\vec{\nu}$ is a finite sequence of ordinals in
 $\mathcal{H}_{\vec{\nu},n}[\Tht](\calk)\cap \vec{\alp}$.
 
By (\ref{eq:controlder}) we have ${\sf k}(\Del)\subset\calh[\Tht]\cap(\calk+1)$ for the set $\Del$ of $\Sig^{1}_{k}(\calk)$-sentences.
Hence for any $\rho\in H_{k}(\vec{\nu},\gam)$,
$\Del^{(\rho,\calk)}$ is a set of $\Sig^{1}_{k}(\rho)$-sentences, and a fortiori of $\Pi^{1}_{k}(\calk)$-sentences with 
 $(\Del^{(\rho,\calk)})^{(\kap,\calk)}=\Del^{(\rho,\calk)}$.

Let
$H_{k}(\vec{\nu},\kap,\hat{a})$ denote the class
{\small
\[
\{\rho\in \bigcap_{i<lh(\vec{\alp})}Mh_{k+i,n}((\vec{\nu}\bullet\vec{\alp})[i] )[\Tht\cup\{\kap\}]\cap\kap
:
\calh_{\hat{a},n}[\Tht\cup\{\kap\}](\rho)\cap\kap\subset\rho 
\}
\]
}

By Proposition \ref{prp:clshull}.\ref{prp:L4.10.1.1a}, $\calh_{\hat{a},n}[\Tht\cup\{\calk\}](\kap)\cap\calk\subset\kap$
 and $\hat{a}\geq \gam$, we have
\beqn\label{eq:L4.10.0}
\rho\in H_{k}(\vec{\nu},\kap,\hat{a}) 
\Rarw \rho\in H_{k}(\vec{\nu},\gam)
\eeqn

Next we have
$\widehat{a_{r}(\rho)}:=\gamma+a_{r}(\rho)\in\mathcal{H}_{\widehat{a_{r}(\rho)},n}[\Tht\cup\{\rho\}](\calk)\cap\hat{a}$ by 
$\widehat{a_{r}(\rho)}\geq\gamma$ and $a_{r}(\rho)<a$.
If $\kap<\calk$,
 Propositions \ref{prp:clshull}.\ref{prp:Mh3} and \ref{prp:clshull}.\ref{prp:L4.10.1.1c} with $a_{r}(\rho)<a$ we have
$\kap\in Mh_{k-1,n}((\widehat{a_{r}(\rho)})*\vec{\alp})[\Tht\cup\{\rho,\kap\}]$, and
$\calh_{\widehat{a_{r}(\rho)},n}[\Tht\cup\{\rho\}](\kap)\subset
\calh_{\hat{a},n}[\Tht](\kap)$ for $\rho<\kap$.

For each $\rho\in H_{k}(\vec{\nu},\kap,\hat{a})$,
IH with (\ref{eq:L4.10.0}) 
yields 

\beqn\label{eq:L4.10case1r}
(\calh_{\hat{a}},\Tht\cup\{\rho,\kap\},{\sf ZFLK}_{k-1,n})\vdash^{\kap+\ome  a_{r}(\rho)}_{\kap}
\Lam,\Gamma^{(\kap,\calk)}, \Del^{(\rho,\calk)}
\eeqn

On the other hand we have 
$\widehat{a_{\ell}}:=\gamma+a_{\ell}\in\mathcal{H}_{\widehat{a_{\ell}},n}[\Tht]\cap\hat{a}$.
Let 
$\rho\in Mh_{k-1,n}((\widehat{a_{\ell}})*\vec{\alp})[\Tht\cup\{\kap\}]\cap H_{k}(\vec{\nu},\kap,\hat{a})$.
Then by Propositions \ref{prp:clshull}.\ref{prp:L4.10.1.1a}  and \ref{prp:clshull}.\ref{prp:L4.10.1.1aa}
we have
$\rho\in Mh_{k-1,n}((\widehat{a_{\ell}})*\vec{\alp})[\Tht\cup\{\rho\}]\cap\kap$,
and
$\calh_{\widehat{a_{\ell}},n}[\Tht](\rho)\cap\calk\subset\rho$.
Hence by IH 
we have for any 
$\rho\in Mh_{k-1,n}((\widehat{a_{\ell}})*\vec{\alp})[\Tht\cup\{\kap\}]\cap H_{k}(\vec{\nu},\kap,\hat{a})$ 
and for any $\del\in\Del$, 
\beqn\label{eq:L4.10case1l}
(\calh_{\hat{a}},\Tht\cup\{\rho\},{\sf ZFLK}_{k-1,n})\vdash^{\rho+\ome  a_{\ell}}_{\rho}\Lam,\Gamma^{(\rho,\calk)}, \lnot\del^{(\rho,\calk)}
\eeqn
Here note that $\rho>\lam$, which is seen as in (\ref{eq:CollapsingthmKR100}).

Now let
\[
M_{\ell} :=
Mh_{k-1,n}((\widehat{a_{\ell}})*\vec{\alp})[\Tht\cup\{\kap\}]\cap H_{k}(\vec{\nu},\kap,\hat{a})
.\]
Then it is clear that
 $M_{\ell}$  is the resolvent class for
$\kap\in Mh_{k-1,n}((\hat{a})*\vec{\alp})[\Tht\cup\{\kap\}]$ 
with respect to $\vec{\mu}=(\widehat{a_{\ell}})*\vec{\nu}$ and $\hat{a}$.

Since $\max\{\rk_{L}(t): t\in{\sf k}(\del^{(\rho,\calk)})\}\leq\rho$, 
we have $\mbox{rk}(\del^{(\rho,\calk)})<\ome(\rho+1)<\kap$ by Proposition \ref{lem:rank}.\ref{lem:rank1}.
From (\ref{eq:L4.10case1r}) and (\ref{eq:L4.10case1l}) by several $(cut)$'s of $\del^{(\rho,\calk)}$ 
 we obtain for $a(\rho)=\max\{a_{\ell},a_{r}(\rho)\}$
and some $p<\ome$
\[
(\calh_{\hat{a}},\Tht\cup\{\kap,\rho\},{\sf ZFLK}_{k-1,n})\vdash^{\kap+\ome  a(\rho)+p}_{\kap}
\Lam,\Gamma^{(\kap,\calk)}, \Gamma^{(\rho,\calk)}
\]
for any $\rho\in M_{\ell}$.

By Inversion lemma \ref{lem:2ndinversion} we obtain for any 
$B\in\calP(\calk)\cap L_{\calk^{+}}$
{\small
\beqn \label{eq:L4.10case1.1a}
\{(\calh_{\hat{a}},\Tht\cup\{\kap,\rho\},{\sf ZFLK}_{k-1,n}))\vdash^{\kap+\ome  a(\rho)+p}_{\kap}
\Lam,\Gamma^{(\kap,\calk)}, \Gamma(R_{B,\calk})^{(\rho,\calk)}: \rho\in M_{\ell}\}
\eeqn
}
where $\Gam(R_{B,\calk})=\{\tau(R_{B,\calk}):(\fal X\subset\calk\,\tau(X))\in\Gam\}$
for $\Sig^{1}_{k-1}(\calk)$-formulae $\tau(X)$,
and $\tau(R_{B,\calk})^{(\rho,\calk)}\equiv \tau^{(\rho,\calk)}(R_{B,\rho})$.

On the other hand
we have by Tautology lemma \ref{lem:tautology}
for each $ \tau(R_{B,\calk})^{(\kap,\calk)}\in\Gam(R_{B,\calk})^{(\kap,\calk)}$
\beqn\label{eq:L4.10case1.1b}
(\calh,\Tht\cup\{\kap\},{\sf ZFLK}_{k-1,n})\vdash^{2\footnotesize{\mbox{{\rm rk}}}(\tau(R_{B,\calk})^{(\kap,\calk)})}_{0} 
\Lam,\Gam(R_{B,\calk})^{(\kap,\calk)}, \lnot\tau(R_{B,\calk})^{(\kap,\calk)}
\eeqn
where 
$2\mbox{{\rm rk}}(\tau(R_{B,\calk})^{(\kap,\calk)}) \leq\kap+p$ for some $p<\ome$ by Proposition \ref{lem:rank}.\ref{lem:rank1}.
Note that $\mbox{{\rm rk}}(\tau(R_{B,\calk})^{(\kap,\calk)})$ is independent from subsets $B$, 
where $\mbox{{\rm rk}}(R_{B,\calk}(c)^{(\kap,\calk)})=\mbox{{\rm rk}}(R_{B,\kap}(c))=\kap$,
and $R_{B,\kap}(c),\Gam$ is an axiom if $c\in B\cap\kap$.

Moreover we have
$\sup\{2\mbox{{\rm rk}}(\tau(R_{B,\calk})^{(\kap,\calk)}), \kap+\ome  a(\rho)+p:\rho\in M_{\ell}\}
\leq \kap+\ome a_{0}+p\in\calh[\Tht\cup\{\kap\}]$
with $\sup\{ a_{\ell}, a_{r}(\rho): \rho\in H_{k}(\vec{\nu},\gam)\}\leq a_{0}\in\calh[\Tht]\cap a$ by (\ref{eq:Mhordgap}).

Since $\lnot\Gam(R_{B,\calk})^{(\kap,\calk)}$ consists of $\Pi^{1}_{k-1}(\kap)$-sentences, 
by an inference rule 
$(\kap\in Mh_{k-1,n}(\hat{a}*\vec{\alp})[\Tht\cup\{\kap\}],\vec{\mu})$
from (\ref{eq:L4.10case1.1b}) and (\ref{eq:L4.10case1.1a}) 
we conclude
\[
(\calh_{\hat{a},n},\Tht\cup\{\kap\},{\sf ZFLK}_{k-1,n})\vdash^{\kap+\ome a_{0}+p+1}_{\kap}\Lam,\Gamma^{(\kap,\calk)}, \Gam(R_{B,\calk})^{(\kap,\calk)}
\]
for any $B\in\calP(\calk)\cap L_{\calk^{+}}$, where
 $\hat{a}\leq\gam_{k-1,n}$ and $\vec{\alp}(i)\leq\gam_{k+i,n}=\vec{\alp}_{k,n}(i)$,
$\kap\in Mh_{k-1,n}(\hat{a}*\vec{\alp})[\Tht\cup\{\kap\}]$ by $\kap\in Mh_{k-1,n}((\hat{a})*\vec{\alp}_{k,n})[\Tht\cup\{\kap\}]$ when
$\kap<\calk$.

Then by several $(\fal^{2}(\kap))$'s we conclude
\[
(\calh_{\hat{a},n},\Tht\cup\{\kap\},{\sf ZFLK}_{k-1,n})\vdash^{\kap+\ome a}_{\kap}\Lam,\Gamma^{(\kap,\calk)}, \Gam^{(\kap,\calk)}
\]
{\bf Case 2}.
Second consider the case when the last inference introduces a $\Pi^{1}_{k}(\calk)$-sentence 
$(\fal X\subset\calk\,\tau(X))\in\Gam$. 
For an $a_{0}<a$ and $J=\calP(\calk)\cap L_{\calk^{+}}$

\[
\infer[(\fal^{2}(\calk))]{(\mathcal{H},\Theta,{\sf ZFLK}_{k,n})\vdash^{a}_{\calk}\Lam,\Gamma}
{
\{
(\mathcal{H},\Theta,{\sf ZFLK}_{k,n})\vdash^{a_{0}}_{\calk}\Lam,\Gamma,
\tau(R_{B,\calk})
:
B\in J)
\}
}
\]
IH yields for each $B\in J$ and $(R_{B,\calk})^{(\kap,\calk)}=R_{B,\kap}$
\[
(\mathcal{H}_{\hat{a}},\Theta\cup\{\kap\},{\sf ZFLK}_{k-1,n})\vdash^{\kap+ \ome a_{0}}_{\kap}\Lam,\Gamma^{(\kap,\calk)},
\tau(R_{B,\calk})^{(\kap,\calk)}
\]
$(\fal^{2}(\kap))$ yields (\ref{eq:CollapsingthmKR01}) with $(\fal X\subset\kap\,\tau(X)^{(\kap,\calk)})\in\Gam^{(\kap,\calk)}$.
\\

\noindent
{\bf Case 3}.
Third consider the case:
for a true literal $M\equiv(R_{B,\sig}(d))$, $M\in\Lam\cup\Gamma$,
where $\sig\leq\calk$ and $d\in B\cap\sig$.
\[
\infer[(\bigwedge)]{(\mathcal{H},\Theta,{\sf ZFLK}_{k,n})\vdash^{a}_{\calk}\Lam,\Gamma}
{}
\]
First consider the case when $\sig<\calk$.
Then $M^{(\kap,\calk)}\equiv M\equiv(R_{B,\sig}(d))\in\Lam\cup\Gamma^{(\kap,\calk)}$.

Second consider the case when $\sig=\calk$.
Then $M^{(\kap,\calk)}\equiv(R_{B,\kap}(d))\in\Gamma^{(\kap,\calk)}$.
It suffices to show $d=\mbox{{\rm rk}}_{L}(d)<\kap$.
We have 
$d\in{\sf k}(R_{B,\calk}(d))\cap\calk\subset\mathcal{H}[\Tht]\cap\calk\subset\kap$
by (\ref{eq:controlder}) and (\ref{eq:CollapsingthmKR100}).
\\

\noindent
{\bf Case 4}.
Fourth consider the case when the last inference introduces a $\Pi^{1}_{0}(\calk)$-sentence 
$(\exi x<\calk\,\tau(x))\in\Gam$.
For a $d<\calk$
\[
\infer[(\bigvee)]{(\calh,\Tht,{\sf ZFLK}_{k,n})\vdash^{a}_{\calk}\Lam,\Gam}
{
(\calh,\Tht,{\sf ZFLK}_{k,n})\vdash^{a_{0}}_{\calk}\Lam,\Gam,\tau(d)
}
\]
where if $(\exi x<\calk\,\tau(x))\in\Sig_{n}$ (this means $(\exi x<\calk\,\tau(x))\in\Del_{0}$), 
then $d=(\mu x<\calk\, \tau(x))\in\calh[\Tht]$,
cf. (\ref{df:assigndc0}) in Definition \ref{eq:dfmu}.

Without loss of generality we can assume that $d\in{\sf k}(\tau(d))\subset\calh[\Tht]$.
Then as in {\bf Case 3} we see that $d<\kap<\kap+ \ome a$, cf. (\ref{eq:bigveebnd}).
Moreover when $(\exi x<\calk\,\tau(x))\in\Sig_{n}$, we have $(\exi x<\calk\,\tau(x))^{(\kap,\calk)}\in\Sig_{n}$ and
$d=(\mu x<\calk\, \tau(x))=(\mu x<\kap\,\tau^{(\kap,\calk)}(x))$
since $L_{\kap}$ is an elementary submodel of $L_{\calk}$, which is seen from $\calh_{\hat{a},n}[\Tht](\kap)\cap\calk\subset\kap$.

IH yields with $(\exi x<\calk\,\tau(x))^{(\kap,\calk)}\equiv\exi x<\kap\, \tau(x)^{(\kap,\calk)}\in\Gam^{(\kap,\calk)}$
\[
\infer[(\bigvee)]{(\calh_{\hat{a}},\Tht\cup\{\kap\},{\sf ZFLK}_{k-1,n})\vdash^{\kap+ \ome a}_{\kap}\Lam,\Gam^{(\kap,\calk)}}
{
(\calh_{\hat{a}},\Tht\cup\{\kap\},{\sf ZFLK}_{k-1,n})\vdash^{\kap+ \ome a_{0}}_{\kap}\Lam,\Gam^{(\kap,\calk)},\tau(d)^{(\kap,\calk)}
}
\]
{\bf Case 5}. 
Fifth consider the case when the last inference introduces a $\Pi^{1}_{0}(\calk)$-sentence 
$(\fal x<\calk\,\tau(x))\in\Gam$.

\[
\infer[(\bigwedge)]{(\mathcal{H},\Theta,{\sf ZFLK}_{k,n})\vdash^{a}_{\calk}\Lam,\Gamma}
{
\{
(\mathcal{H},\Theta\cup\{\alp\},{\sf ZFLK}_{k,n})\vdash^{a(\alp)}_{\calk}\Lam,\Gamma,
\tau(\alp)
:
\alp\in J
\}
}
\]
where if $(\fal x<\calk\,\tau(x))\in\Sig_{n}$, then $J=\{d\}$ for $d=(\mu x<\calk\, \lnot\tau(x))$,
and $J=\calk$ otherwise.
In the former case we see (\ref{eq:CollapsingthmKR01}) from IH and 
$(\mu x<\calk\, \lnot\tau(x))=(\mu x<\kap\,\lnot\tau^{(\kap,\calk)}(x))$.

In what follows suppose $(\fal x<\calk\,\tau(x))\not\in\Sig_{n}$. Then $(\fal x<\kap\,\tau^{(\kap,\calk)}(x))\not\in\Sig_{n}$.
Let $\alp<\kap$.
We have $\calh_{\widehat{a(\alp)},n}[\Tht\cup\{\alp\}](\kap)=\calh_{\widehat{a(\alp)},n}[\Tht](\kap)$,
 and $\kap\in Mh_{k-1,n}((\widehat{a(\alp)})*\vec{\alp}_{k,n}[\Tht\cup\{\alp\}]$ by
$\kap\in Mh_{k-1,n}((\hat{a})*\vec{\alp}_{k,n})[\Tht]$ and Proposition \ref{prp:clshull}.\ref{prp:L4.10.1.1c}
when $\kap<\calk$.

IH yields for $\alp<\kap$ 
\[
(\mathcal{H}_{\hat{a}},\Theta\cup\{\kap,\alp\},{\sf ZFLK}_{k-1,n})\vdash^{\kap+ \ome a(\alp)}_{\kap}\Lam,\Gamma^{(\kap,\calk)},
\tau(\alp)^{(\kap,\calk)}
\]
$(\bigwedge)$ yields (\ref{eq:CollapsingthmKR01}) with
 $\fal x<\kap\,\tau(x)^{(\kap,\calk)}\equiv(\fal x<\calk\,\tau(x))^{(\kap,\calk)}\in\Gam^{(\kap,\calk)}$.
\\

\noindent
{\bf Case 6}.
Sixth consider the case when the last inference introduces a $\Sig^{1}_{k-1}(\calk)$-sentence 
$(\exi X\subset\calk\,\tau(X))\in\Gam$.
For a $B\in\calP(\calk)\cap L_{\calk^{+}}$
\[
\infer[(\exi^{2}(\calk))]{(\calh,\Tht,{\sf ZFLK}_{k,n})\vdash^{a}_{\calk}\Lam,\Gam}
{
(\calh,\Tht,{\sf ZFLK}_{k,n})\vdash^{a_{0}}_{\calk}\Lam,\Gam,\tau(R_{B,\calk})
}
\]
IH with $(R_{B,\calk})^{(\kap,\calk)}\equiv R_{B,\kap}$ yields
\[
\infer[(\exi^{2}(\kap))]{(\calh_{\hat{a}},\Tht\cup\{\kap\},{\sf ZFLK}_{k-1,n})\vdash^{\kap+ \ome a}_{\kap}\Lam,\Gam^{(\kap,\calk)}}
{
(\calh_{\hat{a}},\Tht\cup\{\kap\},{\sf ZFLK}_{k-1,n})\vdash^{\kap+ \ome a_{0}}_{\kap}\Lam,\Gam^{(\kap,\calk)},\tau^{(\kap,\calk)}(R_{B,\kap})
}
\]
{\bf Case 7}.
Seventh consider the case when the last inference is a $(\sig\in Mh_{j,n}(\vec{\alp})[\Tht],\vec{\nu})$
for a $\sig<\calk$ and a $j<N$:

{\small
\[
\infer{(\mathcal{H},\Tht,{\sf ZFLK}_{k,n})\vdash^{a}_{\calk}\Lam,\Gam}
{
\{(\mathcal{H},\Tht,{\sf ZFLK}_{k,n})\vdash^{a_{\ell}}_{\calk}\Lam,\Gamma,\lnot\del\}_{\del\in\Del}
&
 \{(\mathcal{H},\Tht\cup\{\rho\},{\sf ZFLK}_{k,n})\vdash^{a_{r}(\rho)}_{\calk}
 \Lam,\Gamma, \Del^{(\rho,\sig)}\}_{ \rho\in H_{j}(\vec{\nu},\gam)}
}
\]
}
where $\Del$ is a finite set of $\Sig^{1}_{j}(\sig)$-sentences,
and 
$H_{j}(\vec{\nu},\gam)=
\{\rho\in \bigcap_{i<lh(\vec{\alp})}Mh_{j+i,n}((\vec{\nu}\bullet\vec{\alp})[i] [\Tht]\cap\sig : \calh_{\gam,n}[\Tht](\rho)\cap\sig\subset\rho \}$ 
is the resolvent class for $\sig\in Mh_{j,n}(\vec{\alp})[\Tht]$ 
with respect to a $\vec{\nu}$ and $\gam$.

We have $\sig\in\mathcal{H}[\Tht]\cap\calk\subset\kap$
by (\ref{eq:coeff}) and (\ref{eq:CollapsingthmKR100}).
Hence $\sig<\kap$, and
$\Del\subset\Sig^{2}_{0}(\sig)\subset\Del_{0}(\kap)$ and 
$\del^{(\kap,\calk)}\equiv \del$ for any $\del\in\Del$.

Let $H_{j}(\vec{\nu},\kap,\hat{a})
=\{\rho\in \bigcap_{i<lh(\vec{\alp})}Mh_{j+i,n}((\vec{\nu}\bullet\vec{\alp})[i] )[\Tht\cup\{\kap\}]\cap\sig : 
\calh_{\hat{a},n}[\Tht\cup\{\kap\}](\rho)\cap\sig\subset\rho \}$
be a resolvent class for $\sig\in Mh_{j,n}(\vec{\alp})[\Tht\cup\{\kap\}]$.
Then $H_{j}(\vec{\nu},\kap,\hat{a})\subset H_{j}(\vec{\nu},\gam)$ as in {\bf Case 1}.

From IH
we obtain the assertion (\ref{eq:CollapsingthmKR01}) by an inference rule $(\sig\in Mh_{j,n}(\vec{\alp})[\Tht\cup\{\kap\}],\vec{\nu})$
with the resolvent class $H_{j}(\vec{\nu},\kap,\hat{a})$.
\\

\noindent
{\bf Case 8}.
Eighth consider the case when the last inference is a $(cut)$.
\[
\infer[(cut)]{(\calh,\Tht,{\sf ZFLK}_{k,n})\vdash^{a}_{\calk}\Lam,\Gam}
{
 (\calh,\Tht,{\sf ZFLK}_{k,n})\vdash^{a_{0}}_{\calk}\Lam,\Gam,\lnot C
&
 (\calh,\Tht,{\sf ZFLK}_{k,n})\vdash^{a_{0}}_{\calk}C,\Lam,\Gam
}
\]
where $a_{0}<a$ and $\mbox{rk}(C)<\calk$.
Then $C\in\Del_{0}(\calk)$ by Proposition \ref{lem:rank}.\ref{lem:rank3}.
On the other side we have ${\sf k}(C)\subset\calk$
by Proposition \ref{lem:rank}.\ref{lem:rank1}.
Then ${\sf k}(C)\subset\calh[\Tht]\cap\calk\subset\kap$
by (\ref{eq:controlder}) and (\ref{eq:CollapsingthmKR100}).
Hence $C^{(\kap,\calk)}\equiv C$ and
 $\mbox{rk}(C^{(\kap,\calk)})<\kap$ again by Proposition \ref{lem:rank}.\ref{lem:rank1}.

IH yields
\[
 (\calh_{\hat{a}},\Tht\cup\{\kap\},{\sf ZFLK}_{k-1,n})\vdash^{\kap+ \ome a_{0}}_{\kap}\Lam,\Gam^{(\kap,\calk)},\lnot C^{(\kap,\calk)}
 \]
 and
 \[
 (\calh_{\hat{a}},\Tht\cup\{\kap\},{\sf ZFLK}_{k-1,n})\vdash^{\kap+ \ome a_{0}}_{\kap}C^{(\kap,\calk)},\Lam,\Gam^{(\kap,\calk)}
 \]
 Then by a $(cut)$ we obtain
 \[
 (\calh_{\hat{a}},\Tht\cup\{\kap\},{\sf ZFLK}_{k-1,n})\vdash^{\kap+ \ome a}_{\kap}\Lam,\Gam^{(\kap,\calk)}
\]
{\bf Case 9}.
Ninth consider the case when the last inference is an $(\mbox{{\bf F}})$
where
either $F=F^{\Sigma_{1}}_{x\cup\{\lambda\}}$ for a $\lambda\in Reg$
or 
$F=F^{\Sigma_{n}}_{x}$.
Let $F"A_{0}\equiv A\in rng(F)$.
Then $A_{0}\in\Sig_{n}$.
IH yields (\ref{eq:CollapsingthmKR01}).
\\

\noindent
All other cases are seen easily from IH.
\eprf

\bcor\label{cor:onestepdown}{\rm ($\mbox{FiX}^{i}({\sf ZFL})$)}\\
Let $k>0$.
Suppose for the operator $\mathcal{H}_{\gam_{k},n}$ 
\[
(\mathcal{H}_{\gam_{k,n},n},\emptyset,{\sf ZFLK}_{k,n})\vdash^{a_{n}}_{\calk}\Lam,\Gamma
\]
where
$\Lam\subset\Sig^{\Sig_{n}}(\ome_{1})$, 
$\Gamma$ consists of $\Pi^{1}_{k}(\calk)$-sentences.
Then the following holds 

\[
(\mathcal{H}_{\gam_{k-1,n},n},\emptyset,{\sf ZFLK}_{k-1,n})\vdash^{a_{n}}_{\calk}
\Lam,\Gamma
\]
\ecor
\bprf
By Lemma \ref{lem:CollapsingthmKR} with $\gam_{k,n}\in\calh_{\gam_{k,n},n}$ and $\kap=\calk$ we have 
\[
(\mathcal{H}_{\gam_{k,n}+a_{n},n},\{\calk\},{\sf ZFLK}_{k-1,n})\vdash^{\calk+\ome a_{n}}_{\calk}
\Lam,\Gamma
\]
in other words
\[
(\mathcal{H}_{\gam_{k,n}+a_{n},n},\emptyset,{\sf ZFLK}_{k-1,n})\vdash^{\calk+\ome a_{n}}_{\calk}
\Lam,\Gamma
\]
Note that $\gam_{k,n}+a_{n}=\gam_{k-1,n}$ and $\calk+\ome a_{n}=a_{n}$.
\eprf

\bcor\label{cor:onestepdown2}{\rm ($\mbox{FiX}^{i}({\sf ZFL})$)}\\
Let $0\leq k<N$.
Suppose for the operator $\mathcal{H}_{\gam_{N},n}$ 
\[
(\mathcal{H}_{\gam_{N,n},n},\emptyset,\calk^{+},{\sf ZFLK}_{N,n})\vdash^{a_{n}}_{\calk}\Lam,\Gamma
\]
where 
$\Lam\subset\Sig^{\Sig_{n}}(\ome_{1})$, 
$\Gamma$ consists of $\Pi^{1}_{k+1}(\calk)$-sentences.

Then the following holds

\[
(\mathcal{H}_{\gam_{k,n},n},\emptyset,\calk^{+},{\sf ZFLK}_{k,n})\vdash^{a_{n}}_{\calk}
\Lam,\Gamma
\]

\ecor
\bprf
This is seen from Corollary \ref{cor:onestepdown}.
\eprf

\section{Theorems}\label{sec:theorem}

Let us conclude two theorems.

\begin{lemma}\label{th:embedreg}{\rm (Embedding of Axioms)}\\
 For each axiom $A$ in $\mbox{{\rm T}}_{N}(\mathcal{K},I,n)$, there is an $m<\omega$ such that
 for any operator $\mathcal{H}=\calh_{\gam,n}$, the fact that 
 $(\mathcal{H},\emptyset,I,{\sf ZFLK}_{N,n})\vdash^{I\cdot 2}_{I+m}  A$
is provable in $\mbox{{\rm FiX}}^{i}({\sf ZFL})$.

\end{lemma}
\bprf
In this proof let us write $ (\mathcal{H},\Tht)$ for $ (\mathcal{H},\Tht,I,{\sf ZFLK}_{N,n})$.

Let us consider the axiom for $N$-indescribability of $\calk$
\begin{equation}
\renewcommand{\theequation}{\ref{eq:Kord}}
\fal x[x\subset\calk \to \vphi_{0}(x)\to\exi \rho<\calk(
\vphi_{0}^{(\rho,\calk)}(x\cap\rho))]
\end{equation}
\addtocounter{equation}{-1}
where $\vphi_{0}$ is the $\Pi^{1}_{N}(\calk)$-translation of a $\Pi^{1}_{N}$-formula $\vphi$.
Let $\vphi_{1}$ be the $\Pi^{1}_{N}(\calk)$-sentence obtained from $\vphi$.
Suppose that $\vphi$ contains a second-order quantifier.

First let $B\subset\calk$.
Then by Tautology lemma \ref{lem:tautology} with ${\sf k}(\vphi_{1}^{(\rho,\calk)}(R_{B,\rho}))=\{\rho\}$
\[
 (\mathcal{H},\emptyset)\vdash^{2\footnotesize{\mbox{{\rm rk}}}(\vphi_{1}(R_{B,\calk}))}_{0}\lnot\vphi_{1}(R_{B,\calk}),\vphi_{1}(R_{B,\calk})
\]
and for each $\rho<\calk$
\[
 (\mathcal{H},\{\rho\})\vdash^{2\footnotesize{\mbox{{\rm rk}}}(\vphi_{1}^{(\rho,\calk)}(R_{B,\rho}))+1}_{0}\lnot\vphi_{1}^{(\rho,\calk)}(R_{B,\rho}),
 \exi \rho<\calk\, \vphi_{1}^{(\rho,\calk)}(R_{B,\rho})
\]
By Proposition \ref{lem:rank}.\ref{lem:rank1} we have $\rho\leq \rk(\vphi_{1}^{(\rho,\calk)}(R_{B,\rho}))<\ome\rho+\ome$, and hence
$\rk(\vphi_{1}^{(\rho,\calk)}(R_{B,\rho}))<\rk(\vphi_{1}(R_{B,\calk}))<\calk+\ome$ for any $\rho<\calk$.
By the inference rule $(\mathcal{K}\in M_{N})$ we obtain
\[
 (\mathcal{H},\emptyset)\vdash^{\calk+\ome}_{\calk+1}\lnot\vphi_{1}(R_{B,\calk}),
  \exi \rho<\calk\, \vphi_{1}^{(\rho,\calk)}(R_{B,\rho})
\]
On the other hand we have by Proposition \ref{prp:predset} for $d=d(\vphi)<\ome$
\[
(\calh,\{B\})\vdash^{2d}_{0}\lnot\vphi_{0}(B),\vphi_{1}(R_{B,\calk})
\]
and
\[
(\calh,\{B\})\vdash^{2d+2}_{0}\lnot\exi \rho<\calk\, \vphi_{1}^{(\rho,\calk)}(R_{B,\rho}),\exi \rho<\calk\, \vphi_{0}^{(\rho,\calk)}(B\cap\rho)
\]
Two $(cut)$'s yield
\[
(\calh,\{B\})\vdash^{\calk+\ome+2}_{\calk+\ome}\lnot\vphi_{0}(B),\exi \rho<\calk\, \vphi_{0}^{(\rho,\calk)}(B\cap\rho)
\]
and
\[
(\calh,\{B\})\vdash^{\calk+\ome+5}_{\calk+\ome}B\not\subset\calk\lor(\lnot\vphi_{0}(B)\lor \exi \rho<\calk\, \vphi_{0}^{(\rho,\calk)}(B\cap\rho))
\]
Second for $a\not\subset\calk$ we have
$
(\calh,\{a\})\vdash^{1}_{0}a\not\subset\calk
$
and
\[
(\calh,\{a\})\vdash^{2}_{0}a\not\subset\calk\lor(\lnot\vphi_{0}(a)\lor \exi \rho<\calk\, \vphi_{0}^{(\rho,\calk)}(a\cap\rho))
\]
Therefore an inference rule $(\bigwedge)$ yields
\[
(\calh,\emptyset)\vdash^{\calk+\ome+6}_{\calk+\ome}
\fal x[x\subset\calk \to \vphi_{0}(x)\to\exi \rho<\calk(
\vphi_{0}^{(\rho,\calk)}(x\cap\rho))]
\]

Other axioms are seen as in \cite{liftupZF}.
\eprf

\blem\label{th:embedregthm}{\rm (Embedding)}
\benu
\item
If $\mbox{{\rm T}}_{N}(\mathcal{K},I,n)\vdash \Gam$ for sets $\Gam$ of sentences, 
there are $m,k<\ome$ such that for any operator $\calh=\calh_{\gam,n}$, the fact that
 $(\calh,\emptyset,I,{\sf ZFLK}_{N,n})\vdash_{I+m}^{I\cdot 2+k}\Gam$ 
 is provable in $\mbox{{\rm FiX}}^{i}({\sf ZFL})$.
 
 \item
 If ${\sf ZFL}\vdash \Gam$ for sets $\Gam$ of sentences, there are $n,m,k<\ome$ such that
 for any operator $\calh=\calh_{\gam,n}$, the fact that
 $(\calh,\emptyset,I,{\sf ZFLK}_{-2,n})\vdash_{I+m}^{I\cdot 2+k}\Gam$ 
is provable in $\mbox{{\rm FiX}}^{i}({\sf ZFL})$.
 \eenu
\elem

\begin{theorem}\label{thm:1k}
Suppose for a $\Sig^{1}_{k+2}(\calk)$-sentence $\vphi$,
\[
{\sf ZFLK}_{N}\vdash\vphi
\]
Then we can find an $n<\ome$ such that
\[
{\sf ZFLK}_{k,n}\vdash\vphi
\]
In short, ${\sf ZFLK}_{N}$ is $\Sig^{1}_{k+2}(\calk)$-conservative over 
${\sf ZFLK}_{k}=\bigcup_{n<\ome}{\sf ZFLK}_{k,n}={\sf ZFL}+\{\calk\in Mh_{k,n}(\vec{\alp}_{k,n}): n<\ome\}$.
\end{theorem}
\bprf
Suppose a $\Sig^{1}_{k+2}(\calk)$-sentence $\vphi\equiv\exi X\subset\calk\,\tht(X)$ is provable in ${\sf ZFLK}_{N}$.
Let $B:=\mu B\in\calP(\calk)\cap L_{\calk^{+}}(\tht(R_{B,\calk}))$.
By Embedding lemma \ref{th:embedregthm} 
there exist  $m,p<\ome$ such that
$(\calh_{0,n},\emptyset,I,{\sf ZFLK}_{N,n})\vdash^{I\cdot 2+p}_{I+m}\tht(R_{B,\calk})$ with $n=m+3$.
Let us work temporarily in $\mbox{FiX}^{i}({\sf ZFL})$.
By Predicative elimination \ref{lem:predcereg}.\ref{lem:predcereg2} and \ref{lem:predcereg}.\ref{lem:predcereg5}, 
$(\calh_{\ome_{m-1}(I\cdot 2+p),n},\emptyset,I,{\sf ZFLK}_{N,n})\vdash^{\ome_{m}(I\cdot 2+p)}_{I}\tht(R_{B,\calk})$.
By Collapsing lemma \ref{th:Collapsingthmreg1K}, 
$(\calh_{\gam_{N,n},n},\emptyset,\calk^{+},{\sf ZFLK}_{N,n})\vdash^{b_{n}}_{b_{n}}\tht(R_{B,\calk})$
for $b_{n}=\Psi_{\calk^{+},n}(c)$ for $c=\ome_{m+2}(I+1)>\ome_{m+1}(I\cdot 2+p)$.
By Predicative elimination  \ref{lem:predcereg}.\ref{lem:predcereg2} for $a_{n}=\vphi b_{n}b_{n}$,
\[
(\calh_{\gam_{N,n},n},\emptyset,\calk^{+},{\sf ZFLK}_{N,n})\vdash^{a_{n}}_{\calk}\tht(R_{B,\calk})
\]
By Corollary \ref{cor:onestepdown2} 
\[
(\calh_{\gam_{k,n},n},\emptyset,\calk^{+},{\sf ZFLK}_{k,n})\vdash^{a_{n}}_{\calk}\tht(R_{B,\calk})
\]

Note that any inference rule $(\pi\in Mh_{i,n}(\vec{\alp})[\Tht],\vec{\nu})$ is correct for $\pi<\calk$ provably in {\sf ZFL}.
Also 
inference rules $(\calk\in Mh_{k,n}(\vec{\alp})[\Tht],\vec{\nu})$ occur only for $\Tht\subset(\calk+1)$
in the derivation establishing the fact 
$(\calh_{\gam_{k,n},n},\emptyset,\calk^{+},{\sf ZFLK}_{k,n})\vdash^{a_{n}}_{\calk}\tht(R_{B,\calk})$.
Hence $\calk\in Mh_{k,n}(\vec{\alp})[\Tht]$ is equivalent to $\calk\in Mh_{k,n}(\vec{\alp})[\emptyset]$ by
Proposition \ref{prp:clshull}.\ref{prp:L4.10.1.1c}.
Moreover the rule $(\calk\in Mh_{k,n}(\vec{\alp}),\vec{\nu})$ with $\vec{\alp}(i)\leq\vec{\alp}_{k,n}(i)$
 is correct assuming $\calk\in Mh_{k,n}(\vec{\alp}_{k,n})$,
we conclude by induction on $a_{n}<\calk^{+}$ that, 
\[
\mbox{FiX}^{i}({\sf ZFLK}_{k,n})\vdash\tht(B)
\]
and hence by Theorem \ref{th:consvintfix}
\[
{\sf ZFLK}_{k,n}\vdash\vphi
\]
\eprf

\begin{theorem}\label{thm:2}
Suppose for a first-order formula $\vphi$
\[
{\sf ZFLK}_{N}\vdash\exi x\in L_{\ome_{1}}\,\vphi(x)
\]
Then we can find an $n<\ome$ such that
\[
{\sf ZFL}\vdash \exi\kap<\calk(\kap=\Psi_{\calk,n}^{\vec{\alp}_{0,n},\emptyset}(\ome_{n-1}(I+1))) \to
\fal \alp[\alp=\Psi_{\ome_{1},n}(\ome_{n-1}(I+1))\to \exi x\in L_{\alp}\,\vphi(x)]
\]
\end{theorem}
\bprf
As in the proof of Theorem \ref{thm:1k}, we see for $n_{0}=m_{0}+3$ and 
$(\exi x\in L_{\ome_{1}}\,\vphi(x))\in\Sig^{\Sig_{n_{0}}}(\ome_{1}))$
\[
(\calh_{\gam_{1,n_{0}},n_{0}},\emptyset,\calk^{+},{\sf ZFLK}_{1,n_{0}})\vdash^{a_{n_{0}}}_{\calk}
\exi x\in L_{\ome_{1}}\,\vphi(x)
\]
Let $\kap=\Psi_{\calk,n_{0}}^{\vec{\alp}_{0,n_{0}},\emptyset}(\gam_{0,n_{0}})<\calk$.
Then Lemma \ref{lem:CollapsingthmKR} yields for $\gam_{0,n_{0}}=\gam_{1,n_{0}}+a_{n_{0}}$,
$\kap+\ome a_{n_{0}}=a_{n_{0}}$ and $\vec{\alp}_{0,n_{0}}=(\gam_{0,n_{0}})*\vec{\alp}_{1,n_{0}}$ that
\[
(\calh_{\gam_{0,n_{0}},n_{0}},\{\kap\},\calk^{+},{\sf ZFLK}_{0,n_{0}})\vdash^{a_{n_{0}}}_{\kap}
\exi x\in L_{\ome_{1}}\,\vphi(x)
\]
Note that there occurs no inference rule $(\calk\in Mh_{i,n_{0}}(\vec{\alp})[\Tht],\vec{\nu})$ for any $k,\vec{\alp}$
in the derivation
establishing this fact.
In other words we have
\[
(\calh_{\gam_{0,n_{0}},n_{0}},\{\kap\},\calk^{+},{\sf ZFLK}_{-1,n_{0}})\vdash^{a_{n_{0}}}_{\kap}
\exi x\in L_{\ome_{1}}\,\vphi(x)
\]
We see by induction up to $a_{n_{0}}<\calk^{+}$ that
\[
\mbox{FiX}^{i}({\sf ZFL})\vdash  \exi\kap<\calk(\kap=\Psi_{\calk,n_{0}}^{\vec{\alp}_{0,n_{0}},\emptyset}(\gam_{0,n_{0}})) \to 
 \exi x\in L_{\ome_{1}}\,\vphi(x)
\]
Hence by Theorem \ref{th:consvintfix}
\[
{\sf ZFL}\vdash \exi\kap<\calk(\kap=\Psi_{\calk,n_{0}}^{\vec{\alp}_{0,n_{0}},\emptyset}(\gam_{0,n_{0}})) \to 
 \exi x\in L_{\ome_{1}}\,\vphi(x)
\]

Again by Embedding lemma \ref{th:embedregthm}, 
Predicative elimination \ref{lem:predcereg}.\ref{lem:predcereg2} and \ref{lem:predcereg}.\ref{lem:predcereg5}, 
Collapsing lemma \ref{th:Collapsingthmreg1K}, we can find $m$ and $n>n_{0}$ such that
$\exi\kap(\kap=\Psi_{\calk,n_{0}}^{\vec{\alp}_{0,n_{0}},\emptyset}(\gam_{0,n_{0}}))\in\Pi_{n}$ and
\[
(\calh_{b+1,n},\emptyset,\ome_{1},{\sf ZFLK}_{-2,n})\vdash^{\bet}_{\bet}
\exi\kap(\kap=\Psi_{\calk,n_{0}}^{\vec{\alp}_{0,n_{0}},\emptyset}(\gam_{0,n_{0}})) \to \exi x\in L_{\ome_{1}}\,\vphi(x)
\]
for $\bet=\Psi_{\ome_{1},n}(b)$, $b=\ome_{m+1}(I\cdot 2+\ome)$.
Boundedness lemma \ref{lem:boundednessreg} yields
\[
(\calh_{b+1,n},\emptyset,\ome_{1},{\sf ZFLK}_{-2,n})\vdash^{\bet}_{\bet}
\exi\kap(\kap=\Psi_{\calk,n_{0}}^{\vec{\alp}_{0,n_{0}},\emptyset}(\gam_{0,n_{0}})) \to \exi x\in L_{\bet}\,\vphi(x)
\]
Thus
{\small
\[
\mbox{FiX}^{i}({\sf ZFL})\vdash
\exi\kap(\kap=\Psi_{\calk,n_{0}}^{\vec{\alp}_{0,n_{0}},\emptyset}(\gam_{0,n_{0}})) \to
\fal \bet[\bet=\Psi_{\ome_{1},n}(\ome_{m+1}(I\cdot 2+\ome))\to \exi x\in L_{\bet}\,\vphi(x)]
\]
}
and hence by Theorem \ref{th:consvintfix}
\[
{\sf ZFL}\vdash
\exi\kap(\kap=\Psi_{\calk,n_{0}}^{\vec{\alp}_{0,n_{0}},\emptyset}(\gam_{0,n_{0}})) \to
\fal \bet[\bet=\Psi_{\ome_{1},n}(\ome_{m+1}(I\cdot 2+\ome))\to \exi x\in L_{\bet}\,\vphi(x)]
\]
Finally we have for $n>n_{0},m+2$ 
\[
{\sf ZFL}\vdash \Psi_{\ome_{1},n}(\ome_{n-1}(I+1))>\Psi_{\ome_{1},n}(\ome_{m+1}(I\cdot 2+\ome))
\]
and
\beqn\label{eq:thm2}
{\sf ZFL}\vdash \exi\kap<\calk(\kap=\Psi_{\calk,n}^{\vec{\alp}_{0,n},\emptyset}(\gam_{0,n})) 
\to
\exi\kap<\calk(\kap=\Psi_{\calk,n_{0}}^{\vec{\alp}_{0,n_{0}},\emptyset}(\gam_{0,n_{0}})) 
\eeqn
Consider the latter (\ref{eq:thm2}).
Let $n>n_{0}$.
Then $\calh_{\gam,n}[\Tht]\supset\calh_{\gam,n_{0}}[\Tht]$ for any $\gam,\Tht$.
Hence $b_{n}=\Psi_{\calk^{+},n}(\ome_{n-1}(I+1))\geq\Psi_{\calk^{+},n_{0}}(\ome_{n_{0}-1}(I+1))=b_{n_{0}}$,
$\gam_{0,n}>\gam_{0,n_{0}}$, and $\vec{\alp}_{0,n}>\vec{\alp}_{0,n_{0}}$.
Let $\kap=\Psi_{\calk,n}^{\vec{\alp}_{0,n},\emptyset}(\gam_{0,n})<\calk$.
Then 
$\kap\in \bigcap_{i<N}Mh_{i,n}(\vec{\alp}_{i,n})[\{\kap\}]$ and $\calh_{\gam_{0,n},n}(\kap)\cap\calk\subset\kap$.
Hence $\calh_{\gam_{0,n_{0}},n_{0}}(\kap)\cap\calk\subset\kap$.
In general we see by induction on $\kap$ using the definition (\ref{eq:dfMhkh}) that
\[
\kap\in Mh_{i,n}(\vec{\alp})[\Tht] \Rarw
\kap\in Mh_{i,n_{0}}(\vec{\alp})[\Tht] 
\]
Therefore we obtain
$\kap\in \bigcap_{i<N}Mh_{i,n_{0}}(\vec{\alp}_{i,n_{0}})[\{\kap\}]$ by 
Proposition \ref{prp:clshull}.\ref{prp:Mh3} with $\vec{\alp}_{0,n}>\vec{\alp}_{0,n_{0}}$, and
$\Psi_{\calk,n_{0}}^{\vec{\alp}_{0,n_{0}},\emptyset}(\gam_{0,n_{0}}) \leq\Psi_{\calk,n}^{\vec{\alp}_{0,n},\emptyset}(\gam_{0,n})=\kap<\calk$
by the definition (\ref{eq:Psivec}).

Consequently we obtain
{\small
\[
{\sf ZFL}\vdash \exi\kap<\calk(\kap=\Psi_{\calk,n}^{\vec{\alp}_{0,n},\emptyset}(\ome_{n-1}(I+1))) \to
\fal \alp[\alp=\Psi_{\ome_{1},n}(\ome_{n-1}(I+1))\to \exi x\in L_{\alp}\,\vphi(x)]
\]
}
\eprf

\end{document}